# Polynomial over Associative $D$-Algebra

Aleks Kleyn

ABSTRACT. In the paper I considered algebra of polynomials over associative $D$-algebra with unit. Using the tensor notation allows to simplify the representation of polynomial. I considered questions related to divisibility of polynomial of any power over polynomial of power 1.

## CONTENTS

## 1. PREFACE

The theory of polynomials over commutative ring (see the definition of polynomial, for instance, [1], pp. 97 - 98) has statements similar to statements from number theory. Such statements like the theorem on the uniqueness of the decomposition of the polynomial into a product of irreducible polynomials (section [3]-48), the remainder theorem ([1], p. 173, the theorem 1.1) are among these statements.

The theory of polynomials over non-commutative algebra is more difficult ([11], p. 48). In this paper, we attempted to advance a bit in this field.

The possibility to represent a polynomial as

$$p(x) = a_0 + \sum_{k=1}^{n} a_k \circ x^k$$

is an important statement from which a lot of statements of the paper follow. This statement is based on the theorem 4.6 and its corollary 4.7.

Aleks_Kleyn@MailAPS.org.
http://AleksKleyn.dyndns-home.com:4080/    http://arxiv.org/a/kleyn_a_1.
http://sites.google.com/site/AleksKleyn/    http://AleksKleyn.blogspot.com/.

Initially, the concept of a tensor representation of map of free algebra was applied to the notation of a linear map of free algebra over commutative ring.[1] This style of notation allowed make statements of non commutative calculus more clear.

Opportunity to present a polynomial using tensor eliminates the complexity and allows see important properties of polynomial. Using theorems considered in [5], I proved the theorem 6.9. I hope this is first step to study divisibility of polynomials.

Alexandre Laugier was first reader of my paper. I appreciate his helpful comments.

## 2. Conventions

**Convention 2.1.** *I assume sum over index $s$ in expression like*

$$a_{s\cdot 0} x a_{s\cdot 1}$$

$\square$

**Convention 2.2.** *Let $A$ be free algebra with finite or countable basis. Considering expansion of element of algebra $A$ relative basis $\overline{\overline{e}}$ we use the same root letter to denote this element and its coordinates. In expression $a^2$, it is not clear whether this is component of expansion of element a relative basis, or this is operation $a^2 = aa$. To make text clearer we use separate color for index of element of algebra. For instance,*

$$a = a^i e_i$$

$\square$

**Convention 2.3.** *It is very difficult to draw the line between the module and the algebra. Especially since sometimes in the process of constructing, we must first prove that the set $A$ is a module, and then we prove that this set is an algebra. Therefore, to write the element of the module, we will also use the convention 2.2.* $\square$

**Convention 2.4.** *Element of $\Omega$-algebra $A$ is called $A$-number. For instance, complex number is also called $C$-number, and quaternion is called $H$-number.* $\square$

## 3. Zero Divisor of Associative $D$-Algebra

Let $D$ be commutative ring and $A$ be associative $D$-algebra with unit.

**Definition 3.1.** Let $a, b \in A$, $a \neq 0$, $b \neq 0$. If $ab = 0$, then $a$ is called **left zero divisor** and $b$ is called **right zero divisor**.[2] If left zero divisor $a$ is right zero divisor, then $a$ is called **zero divisor**. $\square$

**Theorem 3.2.** *Let $a, b \in A$, $a \neq 0$, $b \neq 0$. The equation $ba = 0$ does not follow from the equation $ab = 0$.*[3]

*Proof.* Let $A$ be algebra of $3 \times 3$ matrices. Let

$$E_{12} = \begin{pmatrix} 0 & 1 & 0 \\ 0 & 0 & 0 \\ 0 & 0 & 0 \end{pmatrix} \qquad\qquad E_{23} = \begin{pmatrix} 0 & 0 & 0 \\ 0 & 0 & 1 \\ 0 & 0 & 0 \end{pmatrix}$$

---

[1]See, for instance, section [6]-1, [10]-1.

[2]See also the definition [2]-10.17.

[3]The proof of the theorem is based on the remark in [2] after the definition 10.17.

It is evident that

$$E_{23}E_{12} = \begin{pmatrix} 0 & 0 & 0 \\ 0 & 0 & 0 \\ 0 & 0 & 0 \end{pmatrix}$$

However

$$E_{12}E_{23} = \begin{pmatrix} 0 & 0 & 1 \\ 0 & 0 & 0 \\ 0 & 0 & 0 \end{pmatrix}$$

It is easy to see that both matrices $(E_{12}, E_{23})$ are zero divisors. $\square$

**Theorem 3.3.** *Let $a$ be a right zero divisor of $D$-algebra $A$. Let $b$ be left zero divisor of $D$-algebra $A$. Let $ab \neq 0$. Then for any $d \in A$, $adb$ is zero divisor of $D$-algebra $A$.*

*Proof.* Since $a$ is right zero divisor of $D$-algebra $A$, then there exists $c \neq 0$ such that $ca = 0$. Then

$$c(adb) = (ca)(db) = 0(db) = 0$$

Therefore, $adb$ is right zero divisor of $D$-algebra $A$.

Since $b$ is left zero divisor of $D$-algebra $A$, then there exists $c \neq 0$ such that $bc = 0$. Then

$$(adb)c = (ad)(bc) = (ad)0 = 0$$

Therefore, $adb$ is left zero divisor of $D$-algebra $A$. $\square$

**Theorem 3.4.** *There exists $D$-algebra where left zero divisor is not right zero divisor.*[4]

*Proof.* Let $A$ be free $R$-vector space which has Hamel basis $\overline{\overline{e}}$.[5] Consider $R$-algebra of linear maps[6] $\mathcal{L}(R; A; A)$. Let map $f \in \mathcal{L}(R; A; A)$ be defined by the equation

(3.1)
$$\begin{cases} f \circ e_i = e_{i-1} & i > 1 \\ f \circ e_1 = 0 \end{cases}$$

Let map $g \in \mathcal{L}(R; A; A)$ be defined by the equation

(3.2)
$$g \circ e_i = e_{i+1}$$

Let map $p \in \mathcal{L}(R; A; A)$ be defined by the equation

(3.3)
$$\begin{cases} p \circ e_i = 0 & i > 1 \\ p \circ e_1 = e_1 \end{cases}$$

---

[4]See also example [2]-10.16.

[5]Hamel basis was considered in the definition [9]-2.3.1.

[6]Let $f, g \in \mathcal{L}(R; A; A)$. The sum of maps $f$ and $g$ is defined by the equation

$$(f + g) \circ x = f \circ x + g \circ x$$

The product of maps $f$ and $g$ is defined by the equation

$$(f \circ g) \circ x = f \circ (g \circ x)$$

From equations (3.1), (3.2), it follows that[7]

(3.4)                                    $f \circ g = 1$

3.4.1: From equations (3.1), (3.3), it follows that

$$(f \circ p) \circ e_1 = f \circ (p \circ e_1) = f \circ e_1 = 0$$

$$(f \circ p) \circ e_i = f \circ (p \circ e_i) = f \circ 0 = 0 \qquad \qquad i > 1$$

Therefore, the map $f$ is left zero divisor.

3.4.2: Let $h \in \mathcal{L}(R; A; A)$ be such map that

(3.5)                                    $h \circ f = 0$

From equations (3.4), (3.5), it follows that

$$0 = 0 \circ g = (h \circ f) \circ g = h \circ (f \circ g) = h \circ 1 = h$$

Therefore, the map $f$ is not right zero divisor.

3.4.3: From equations (3.2), (3.3), it follows that

$$(p \circ g) \circ e_i = p \circ (g \circ e_i) = p \circ e_{i+1} = 0$$

Therefore, the map $g$ is right zero divisor.

3.4.4: Let $h \in \mathcal{L}(R; A; A)$ be such map that

(3.6)                                    $g \circ h = 0$

From equations (3.4), (3.6), it follows that

$$0 = f \circ 0 = f \circ (g \circ h) = (f \circ g) \circ h = 1 \circ h = h$$

Therefore, the map $g$ is not left zero divisor.

$\square$

**Theorem 3.5.** *There exists D-algebra where left zero divisor is invertible from right.*[8]

---

[7]In the equation (3.4) we see how properties of matrix change when we consider a matrix with countable set of rows and columns instead of matrix with finite number of rows and columns. Relative to the basis $\overline{\overline{e}}$, the map $f$ has matrix

$$f = \begin{pmatrix} 0 & 1 & 0 & 0 & ... \\ 0 & 0 & 1 & 0 & ... \\ 0 & 0 & 0 & 1 & ... \\ ... & ... & ... & ... & ... \end{pmatrix}$$

First column of matrix $f$ linearly depends from other columns of matrix. However, the matrix $f$ is invertible from right.

It is possible that this phenomenon is associated with the statement that a countable set has proper countable subset.

[8]Let $A$ be D-algebra and $a \in A$ be left zero divisor. Then there exists $b \in A$, $b \neq 0$, such that $ab = 0$. Let there exists $c \in A$ such that $ca = 1$. Therefore,

$$b = 1b = (ca)b = c(ab) = c0 = 0$$

From this contradiction, it follows that $a$ is not invertible from left.

*Proof.* Let $\mathcal{L}(R; A; A)$ be $R$-algebra considered in the proof of the theorem 3.4. Let map $f \in \mathcal{L}(R; A; A)$ be defined by the equation (3.1). Let map $g \in \mathcal{L}(R; A; A)$ be defined by the equation (3.2). According to the statement 3.4.1, the map $f$ is left zero divisor. According to the equation (3.4), the map $f$ is invertible from right. □

**Theorem 3.6.** *Let $a$ be a left zero divisor of $D$-algebra $A$. Non zero element of right ideal $Aa$ is left zero divisor of $D$-algebra $A$.*

*Proof.* According to the definition 3.1, there exists $b \in A$, $b \neq 0$, such that $ab = 0$. Then for any $c \in A$

$$(ca)b = c(ab) = c0 = 0$$

Therefore, if $ca \neq 0$, then $ca$ is left zero divisor. □

**Theorem 3.7.** *If we can represent left zero divisor $a$ of $D$-algebra $A$ as product $a = cd$, then either $c$, or $d$ is a left zero divisor.*

*Proof.* If $d$ is left zero divisor, then, according to the theorem 3.6, $a$ is left zero divisor. So to prove the theorem, we consider the case when $d$ is not left zero divisor. According to the definition 3.1, there exists $b \in A$, $b \neq 0$, such that $ab = 0$. Then

$$0 = ab = (cd)b = c(db)$$

According to the definition 3.1, $db \neq 0$. Therefore, $c$ is left zero divisor. □

**Theorem 3.8.** *Let neither $a \in A$, nor $b \in A$ be left zero divisors of $D$-algebra $A$. Then their product $ab$ is not left zero divisor of $D$-algebra $A$.*

*Proof.* If the product $ab$ is a left zero divisor of $D$-algebra $A$, then, according to the theorem 3.7, then either $a$, or $b$ is a left zero divisor. This contradiction proves the theorem. □

To see better the structure of the set of zero divisors of $D$-algebra $A$, we consider the following theorem.

**Theorem 3.9.** *Let $A$ be finite dimensional $D$-algebra. Let $\overline{\overline{e}}$ be a basis of $D$-algebra $A$. Let $C_{kl}^i$ be structural constants of $D$-algebra $A$ relative to the basis $\overline{\overline{e}}$. Then coordinates $a^i$ of left zero divisor*

$$(3.7) \qquad a = a^i e_i$$

*satisfy to equation*

$$(3.8) \qquad \det \|C_{kl}^i a^k\| = 0$$

*Proof.* Since $A$-number $a \neq 0$ is left zero divisor, then according to the definition 3.1, $a \neq 0$ and there exists $A$-number $b \neq 0$

$$b = b^i e_i$$

such that $ab = 0$. Therefore, coordinates of $A$-numbers $a$ and $b$ satisfy to the system of equations

$$(3.9) \qquad C_{kl}^i a^k b^l = 0$$

If we assume that we know $a$, then we can consider the system of equations (3.9) as system of linear equations relative to coordinates $b^l$. The equation (3.8) follows from the statement that number of equations in the system of linear equations

(3.9) equals to the number of unknown and the system of linear equations (3.9) has nontrivial solution. □

## 4. Polynomial over Associative $D$-Algebra

Let $D$ be commutative ring and $A$ be associative $D$-algebra with unit.

**Theorem 4.1.** *Let* $p_k(x)$ *be* **monomial of power** $k$ *over $D$-algebra $A$.*[9] *Then*

4.1.1: *Monomial of power* 0 *has form* $p_0(x) = a_0$, $a_0 \in A$.

4.1.2: *If* $k > 0$, *then*

$$p_k(x) = p_{k-1}(x)xa_k$$

*where* $a_k \in A$.

*Proof.* We prove the theorem by induction over power $n$ of monomial.

Let $n = 0$. We get the statement 4.1.1 since monomial $p_0(x)$ is constant.

Let $n = k$. Last factor of monomial $p_k(x)$ is either $a_k \in A$, or has form $x^l$, $l \geq 1$. In the later case we assume $a_k = 1$. Factor preceding $a_k$ has form $x^l$, $l \geq 1$. We can represent this factor as $x^{l-1}x$. Therefore, we proved the statement. □

**Definition 4.2.** We denote $A_k[x]$ Abelian group generated by the set of monomials of power $k$. Element $p_k(x)$ of Abelian group $A_k[x]$ is called **homogeneous polynomial of power** $k$. □

*Remark* 4.3. According to the definition 4.2, homogeneous polynomial $p_k(x)$ of power $k$ is sum of monomials of power $k$

$$p_k(x) = \sum_s p_{k \cdot s}(x)$$

Let

$$q_k(x) = \sum_t q_{k \cdot t}(x)$$

be homogeneous polynomials of power $k$. Since sum in Abelian group $A_k[x]$ is commutative and associative, then sum of homogeneous polynomials $p$ and $q$ is sum of monomials of power $k$

$$(p + q)(x) = p(x) + q(x) = \sum_s p_{k \cdot s}(x) + \sum_t q_{k \cdot t}(x)$$

□

**Theorem 4.4.** *Abelian group $A_k[x]$ is $A$-module.*

*Proof.* Let

$$p(x) = p_0 x p_1 ... p_{k-1} x p_k$$

be monomial. For any tensor $a \otimes b \in A \otimes A$, there exists transformation of monomial

(4.1)                    $(a \otimes b) \circ p(x) = ap(x)b = ap_0 x p_1 ... p_{k-1} x p_k b$

The set of transformations (4.1) generates representation of $D$-algebra $A \otimes A$ in Abelian group $A_k[x]$. Therefore, Abelian group $A_k[x]$ is $A$-module. □

---

[9]You can see similar definition of monomial over division ring in the section [10]-16. You can see similar definition of monomial over Banach algebra in the section [8]-5.2.

The set of monomials of power $k$ is not a basis of $A$-module $A_k[x]$. For instance,

$$adxb = axdb \quad d \in D$$

For polynomial, we will use notation

$$(a_0, a_1, ..., a_k) \circ x^k = a_0 x a_1 x ... x a_k$$

**Theorem 4.5.** *The map*

$$f : A^{k+1} \to A_k[x]$$

*defined by the equation*

$$f(a_0, a_1, ..., a_k) = (a_0, a_1, ..., a_k) \circ x^k$$

*is polylinear map.*

*Proof.* Let $d \in D$. From equations

$$f(a_0, ..., da_i, ..., a_k) = a_0 x ... (da_i) ... x a_k = d(a_0 x ... a_i ... x a_k)$$
$$= df(a_0, ..., a_i, ..., a_k)$$
$$f(a_0, ..., a_i + b_i, ..., a_k) = a_0 x ... (a_i + b_i) ... x a_k$$
$$= (a_0 x ... a_i ... x a_k) + (a_0 x ... b_i ... x a_k)$$
$$= f(a_0, ..., a_i, ..., a_k) + f(a_0, ..., b_i, ..., a_k)$$

it follows that the map $f$ is linear with respect to $a_i$. Therefore, the map $f$ is polylinear map. $\square$

**Theorem 4.6.** *There exists linear map*

$$(a_0 \otimes a_1 \otimes ... \otimes a_k) \circ x^k = (a_0, a_1, ..., a_k) \circ x^k$$

*Proof.* The theorem follows from theorems [7]-3.6.4, 4.5. $\square$

**Corollary 4.7.** *We can present homogeneous polynomial $p(x)$ in the following form*

$$p(x) = a_k \circ x^k \quad a_k \in A^{\otimes(k+1)}$$

$\square$

**Definition 4.8.** We denote

$$A[x] = \bigoplus_{n=0}^{\infty} A_n[x]$$

direct sum[10] of $A$-modules $A_n[x]$. An element $p(x)$ of $A$-module $A[x]$ is called **polynomial** over $D$-algebra $A$. $\square$

Therefore, we can present polynomial of power $n$ in the following form

$$p(x) = a_0 + a_1 \circ x + ... + a_n \circ x^n \quad a_i \in A^{\otimes(i+1)} \quad i = 0, ..., n$$

**Definition 4.9.** Let

$$p(x) = a_0 + a_1 \circ x + ... + a_n \circ x^n \quad a_i \in A^{\otimes(i+1)} \quad i = 0, ..., n$$

be polynomial of power $n$ over $D$-algebra $A$. $A^{\otimes(i+1)}$-number $a_i$ is called **coefficient of polynomial** $p(x)$. $A^{\otimes(n+1)}$-number $a_n$ is called **leading coefficient of polynomial** $p(x)$. $\square$

---

[10]See the definition of direct sum of modules in [1], page 128. On the same page, Lang proves the existence of direct sum of modules.

## 5. Operations over Polynomials

**Definition 5.1.** Let

$$p(x) = p_0 + p_1 \circ x + ... + p_n \circ x^n$$
$$q(x) = q_0 + q_1 \circ x + ... + q_n \circ x^n$$

be polynomials of power $n$.[11] **Sum of polynomials** is defined by equation

$$(p + q)(x) = p_0 + q_0 + (p_1 + q_1) \circ x + ... + (p_n + q_n) \circ x^n$$

$\square$

**Definition 5.2.** Bilinear map

$$\underline{\otimes} : A^{\otimes n} \times A^{\otimes m} \to A^{\otimes(n+m-1)}$$

is defined by the equation

$$(a_1 \otimes ... \otimes a_n) \underline{\otimes} (b_1 \otimes ... \otimes b_n) = a_1 \otimes ... \otimes a_{n-1} \otimes a_n b_1 \otimes b_2 \otimes ... \otimes b_n$$

$\square$

**Theorem 5.3.** *Homogeneous polynomial $p \circ x^n$, $p \in A^{\otimes(n+1)}$, generates the linear map*

$$A_m[x] \to A_{n+m}[x]$$

*defined by the equation*

(5.1)
$$(p \circ x^n) \circ (q \circ x^m) = (p \underline{\otimes} q) \circ x^{n+m}$$

*Proof.* According to the corollary 4.7, $p \in A^{\otimes(n+1)}$, $q \in A^{\otimes(m+1)}$. According to the definition 5.2,

(5.2)
$$p \underline{\otimes} q \in A^{\otimes((n+1)+(m+1)-1)} = A^{\otimes((n+m)+1)}$$

According to the corollary 4.7 and the equation (5.2), the map (5.1) is defined properly. According to the definition 5.2, from equations

$$\begin{aligned}
(p \circ x^n) \circ (r \circ x^m + s \circ x^m) &= (p \circ x^n) \circ ((r + s) \circ x^m) \\
&= (p \underline{\otimes}(r + s)) \circ x^{n+m} \\
&= (p \underline{\otimes} r) \circ x^{n+m} + (p \underline{\otimes} s) \circ x^{n+m} \\
&= (p \circ x^n) \circ (r \circ x^m) + (p \circ x^n) \circ (s \circ x^m) \\
(p \circ x^n) \circ (d(r \circ x^m)) &= (p \circ x^n) \circ ((dr) \circ x^m) \\
&= (p \underline{\otimes}(dr)) \circ x^{n+m} = d((p \underline{\otimes} r) \circ x^{n+m}) \\
&= d((p \circ x^n) \circ (r \circ x^m))
\end{aligned}$$

where $d \in D$, it is follows that the map (5.1) is linear map.                     $\square$

---

[11]Let $q(x)$ be polynomial of power $m$, $m < n$,

$$q(x) = q_0 + q_1 \circ x + ... + q_m \circ x^m$$

If we assume $q_{m+1} = 0, ..., q_n = 0$, then we can consider the polynomial $q(x)$ as polynomial of power $n$.

*Remark* 5.4. Let

$$p = p_0 \otimes p_1 \otimes ... \otimes p_n$$

$$r = r_0 \otimes r_1 \otimes ... \otimes r_m$$

Then we can write the equation (5.1) in following form

$$
\begin{aligned}
(p_0 x ... x p_n) \circ (r_0 x ... x r_m) &= ((p_0 \otimes ... \otimes p_n) \circ x^n) \circ ((r_0 \otimes ... \otimes r_m) \circ x^m) \\
&= ((p_0 \otimes ... \otimes p_n) \underline{\otimes} (r_0 \otimes ... \otimes r_m)) \circ x^{n+m} \\
&= (p_0 \otimes ... \otimes p_n r_0 \otimes ... \otimes r_m) \circ x^{n+m} \\
&= p_0 x ... x p_n r_0 x ... x r_m
\end{aligned}
$$

Therefore, the equation (5.1) is the definition of product of homogeneous polynomials. □

**Definition 5.5.** Product of homogeneous polynomials $p \circ x^n$, $r \circ x^m$ is defined by the equation

$$(p \circ x^n)(r \circ x^m) = (p \underline{\otimes} r) \circ x^{n+m}$$

□

**Theorem 5.6.** *Let*

$$p(x) = a_k \circ x^k \quad a_k \in A^{\otimes(k+1)}$$

*be monomial of power $k > 1$. Then the polynomial $p(x)$ can be represented using one of the following forms*

(5.3) $$p(x) = (a_{k \cdot 0} \circ x^{k-1})((1 \otimes a_{k \cdot 1}) \circ x)$$

(5.4) $$p(x) = ((a_{k \cdot 0} \circ x^{k-1}) \otimes a_{k \cdot 1}) \circ x$$

*where*

(5.5) $$a_k = a_{k \cdot 0} \underline{\otimes} (1 \otimes a_{k \cdot 1}) \quad a_{k \cdot 0} \in A^{\otimes k} \quad a_{k \cdot 1} \in A$$

*Proof.* Based on the statement 4.1.2 and the theorem 4.6, we can write the monomial $p(x)$ as

(5.6) $$p_k(x) = (a_{k \cdot 0} \circ x^{k-1}) x a_{k \cdot 1}$$

The equation (5.5) follows from the definition 5.2 and the equation

$$a_k = a_{k \cdot 0} \otimes a_{k \cdot 1}$$

Since for given value of $x$, the expression $a_{k \cdot 0} \circ x^{k-1}$ is $A$-number, then the equation (5.4) follows from the equation (5.6) and the theorem 4.6. The equation (5.3) follows from equations (5.5), (5.6) and the definition 5.5. □

**Theorem 5.7.** *Let*

$$p(x) = a_k \circ x^k \quad a_k \in A^{\otimes(k+1)}$$

*be homogeneous polynomial of power $k > 1$. Then the polynomial $p(x)$ can be represented using one of the following forms*

(5.7) $$p(x) = (a_{k \cdot 0 \cdot s} \circ x^{k-1})((1 \otimes a_{k \cdot 1 \cdot s}) \circ x)$$

(5.8) $$p(x) = ((a_{k \cdot 0 \cdot s} \circ x^{k-1}) \otimes a_{k \cdot 1 \cdot s}) \circ x$$

*where*

(5.9)        $a_k = a_{k \cdot 0 \cdot s} \underline{\otimes} (1 \otimes a_{k \cdot 1 \cdot s})$     $a_{k \cdot 0 \cdot s} \in A^{\otimes k}$     $a_{k \cdot 1 \cdot s} \in A$

*Proof.* According to the remark 4.3, homogeneous polynomial of power $k$ is sum of monomials of power $k$. The theorem follows from the theorem 5.6, if we consider induction over number of terms.                                                    $\square$

*Remark* 5.8. In this section, it is unimportant for us whether a polynomial $p(x)$ is zero divisor. We consider in the section 6 the question about zero divisors of $A$-algebra $A[x]$. In the proof of the theorems in this section, we, without loss of generality, assume that the considered polynomials are not zero divisors.     $\square$

Since product of homogeneous polynomials is bilinear map, then the definition 5.5 can be extended to product of any polynomials.

**Theorem 5.9.** *Let*

$$p(x) = p_0 + p_1 \circ x + ... + p_n \circ x^n$$
$$q(x) = q_0 + q_1 \circ x + ... + q_m \circ x^m$$

*be polynomials. If polynomial*

$$r(x) = r_0 + r_1 \circ x + ... + r_k \circ x^k$$

*is* **product of polynomials**

(5.10)                          $r(x) = p(x)q(x)$

*then*

$$k = n + m$$

(5.11)       $r_h = \sum_{i+j=h} p_i \underline{\otimes} q_j$     $h = 0, ..., k$   $i \le n$   $j \le m$

*Proof.* We prove the theorem by induction over $n$, $m$.

- Following lemma follows from the definition 5.5.

   *Lemma* 5.10. The theorem 5.9 is true for product of homogeneous polynomials.

- Let the theorem be true for homogeneous polynomial $p(x) = p_n \circ x^n$ of power $n$ and polynomial $q(x)$ of power $m = l$. Polynomial

   $$q(x) = q_0 + q_1 \circ x + ... + q_l \circ x^l + q_{l+1} \circ x^{l+1}$$

   of power $l + 1$ can be written in the following form

   $$q(x) = q'(x) + q_{l+1} \circ x^{l+1}$$

   where

   $$q'(x) = q_0 + q_1 \circ x + ... + q_l \circ x^l$$

   is polynomial of power $l$. Since product of polynomials is bilinear map, then

   (5.12)       $r(x) = (p_n \circ x^n)q(x) = (p_n \circ x^n)(q'(x) + q_{l+1} \circ x^{l+1})$
   $$= (p_n \circ x^n)q'(x) + (p_n \circ x^n)(q_{l+1} \circ x^{l+1})$$

According to the induction assumption, the power of polynomial $(p_n \circ x^n)q'(x)$ equals $n + l$, as well

$$(5.13) \qquad r_h = \sum_{i+j=h} p_i \underline{\otimes} q_j \quad i = n \quad h = 0, ..., n + l$$

According to the definition 5.5, the power of polynomial $(p_n \circ x^n)(q_{l+1} \circ x^{l+1})$ equals

$$k = n + (l + 1)$$

as well

$$(5.14) \qquad r_h = \sum_{i+j=h} p_i \underline{\otimes} q_j \quad i = n \quad h = n + l + 1$$

Therefore, the theorem is true for homogeneous polynomial $p(x) = p_n \circ x^n$ of power $n$ and polynomial $q(x)$ of power $m = l + 1$. Therefore, we proved the following lemma.

*Lemma* 5.11. The theorem 5.9 is true for product of homogeneous polynomial $p(x)$ and polynomial $q(x)$.

Let the theorem be true for polynomial $p(x)$ of power $n = l$ and polynomial $q(x)$ of power $m$. Polynomial

$$p(x) = p_0 + p_1 \circ x + ... + p_l \circ x^l + p_{l+1} \circ x^{l+1}$$

of power $l + 1$ can be written in the following form

$$p(x) = p'(x) + p_{l+1} \circ x^{l+1}$$

where

$$p'(x) = p_0 + p_1 \circ x + ... + p_l \circ x^l$$

is polynomial of power $l$. Since product of polynomials is bilinear map, then

$$(5.15) \qquad \begin{aligned} r(x) = p(x)q(x) &= (p'(x) + p_{l+1} \circ x^{l+1})q(x) \\ &= p'(x)q(x) + (p_{l+1} \circ x^{l+1})q(x) \end{aligned}$$

According to the induction assumption, the power of polynomial

$$r'(x) = p'(x)q(x)$$

equals $l + m$, as well

$$(5.16) \qquad r'_h = \sum_{i+j=h} p_i \underline{\otimes} q_j \quad i < l + 1 \quad h = 0, ..., l + m$$

According to the lemma 5.11, the power of polynomial

$$r''(x) = (p_{l+1} \circ x^{l+1})q(x)$$

equals

$$k = (l + 1) + m$$

as well

$$(5.17) \qquad r''_h = \sum_{i+j=h} p_i \underline{\otimes} q_j \quad i = l + 1 \quad h = 0, ..., l + m + 1$$

Since

$$r(x) = r'(x) + r''(x)$$

then from the definition 5.1 and from equations (5.16), (5.17), it follows that

$$(5.18) \qquad r_h = \sum_{i+j=h} p_i \underline{\otimes} q_j \quad i \le l+1 \quad h = 0, ..., l+m+1$$

Therefore, the theorem is true for polynomial $p(x)$ of power $n = l+1$ and polynomial $q(x)$ of power $m$. $\qquad \square$

**Theorem 5.12.** *$A$-module $A[x]$ equipped by the product (5.10) is $A$-algebra which is called $A$-algebra of polynomials over $D$-algebra $A$.*

*Proof.* The theorem follows from definitions [7]-3.2.1, 5.5 and the theorem 5.9. $\quad \square$

## 6. Division of Polynomials

In this section, we assume that the algebra $A$ is defined over field $F$ and the algebra $A$ is finite dimensional $F$-algebra.[12] Let $\overline{\overline{e}}$ be the basis of algebra $A$ over field $F$ and $C_{ij}^k$ be structural constants of algebra $A$ relative to the basis $\overline{\overline{e}}$.

Consider the linear equation[13]

$$(6.1) \qquad a \circ x = b$$

where $a = a_{s \cdot 0} \otimes a_{s \cdot 1} \in A^{\otimes 2}$. According to the theorem [4]-10.1.9, we can write the equation (6.1) in standard form

$$(6.2) \qquad a^{ij} e_i x e_j = b$$

$$a^{ij} = a_{s \cdot 0}{}^{i} a_{s \cdot 1}{}^{j} \quad a_{s \cdot 0} = a_{s \cdot 0}{}^{i} e_i \quad a_{s \cdot 1} = a_{s \cdot 1}{}^{i} e_i$$

According to the theorem [4]-10.1.10 equation (6.2) is equivalent to equation

$$(6.3) \qquad a_i^j x^i = b^j$$

$$a_i^j = a^{kr} C_{ki}^p C_{pr}^j \quad x = x^i e_i \quad b = b^i e_i$$

According to the theory of linear equations over field, if determinant

$$(6.4) \qquad \det \| a_i^j \| \ne 0$$

then equation (6.1) has only one solution.

**Definition 6.1.** The tensor $a \in A^{\otimes 2}$ is called **nonsingular tensor** if this tensor satisfies to condition (6.4). $\qquad \square$

**Theorem 6.2.** *Let $F$ be a field. Let $A$ be finite dimensional $F$-algebra. For the linear equation*

$$(6.5) \qquad a \circ x = 0$$

*where $a = a_{s \cdot 0} \otimes a_{s \cdot 1} \in A^{\otimes 2}$. , any $x \in A$ is root iff*

$$(6.6) \qquad a_{s \cdot 0}{}^{k} a_{s \cdot 1}{}^{r} C_{ki}^p C_{pr}^j = 0$$

---

[12]From the proof of the theorem 3.4, we see that statements of linear algebra in a vector space with countable basis are different from similar statements in finite dimensional vector space. Additional research is required before we can state theorems of the section 6 for $F$-algebra with countable basis.

[13]I consider the solving of the equation (6.1) the same way as I have done in the section [5]-4.

*Proof.* According to the theorem [4]-10.1.9, we can write the equation (6.5) in standard form

$$(6.7) \qquad a^{ij} e_i x e_j = 0$$

$$a^{ij} = a_{s \cdot 0}{}^{i} a_{s \cdot 1}{}^{j} \quad a_{s \cdot 0} = a_{s \cdot 0}{}^{i} e_i \quad a_{s \cdot 1} = a_{s \cdot 1}{}^{i} e_i$$

According to the theorem [4]-10.1.10 equation (6.7) is equivalent to equation

$$(6.8) \qquad a_i^j x^i = 0$$

$$a_i^j = a^{kr} C_{ki}^p C_{pr}^j \quad x = x^i e_i$$

Any $x \in A$ is root of the system of linear equations (6.8), iff

$$(6.9) \qquad a_i^j = 0$$

From equations (6.8), (6.9), it follows that

$$(6.10) \qquad a^{kr} C_{ki}^p C_{pr}^j = 0$$

The equation (6.6) follows from equations (6.7), (6.10). $\qquad \square$

It is difficult to say whether the condition (6.6) to be true in some algebra. However we can study particular case of the theorem 6.2.

**Theorem 6.3.** *Let $F$ be a field. Let $A$ be finite dimensional $F$-algebra. In order for any $x \in A$ to be the root of the equation*

$$axb = 0$$

*it is necessary that $a$ is left divisor of $F$-algebra $A$.*

*Proof.* Let in the equation (6.5) $s = 1$, $a_{1 \cdot 0} = a$, $a_{1 \cdot 1} = b$. Let $\overline{\overline{e}}$ be the basis of $F$-algebra $A$. Then the equation (6.6) gets form

$$(6.11) \qquad a^k b^r C_{ki}^p C_{pr}^j = 0$$

If we assume that we know $a$, then we can consider the system of equations (6.11) as system of linear equations relative to coordinates $b^r$. Number of equations in the system of linear equations (6.11) equals to the number of unknown. Since the system of linear equations (6.11) has nontrivial solution, then

$$(6.12) \qquad \det \| a^k C_{ki}^p C_{pr}^j \| = 0$$

for any $r$. From the equation (6.12) it follows that

$$(6.13) \qquad \det \| C_{pr}^j \| \det \| a^k C_{ki}^p \| = 0$$

From the equation (6.13) it follows that either

$$(6.14) \qquad \det \| C_{pr}^j \| = 0$$

or

$$(6.15) \qquad \det \| a^k C_{ki}^p \| = 0$$

Let equation (6.14) be true. Let $e_1 \in \overline{\overline{e}}$, $e_1 = 1$. Then

$$(6.16) \qquad e_p = e_p e_1 = C_{p1}^j e_j$$

From the equation (6.14) it follows that there exist $c^p$, $c \neq 0$, such that

$$(6.17) \qquad c^p C_{p1}^j = 0$$

From equations (6.16), (6.17), it follows that

$$(6.18) \qquad c^p e_p = c^p C_{p1}^j e_j = 0$$

Therefore, vectors $e_k$ are linear dependent. From contradiction, it follows that the equation (6.14) is not true in $F$-algebra $A$ with unit.

From the equation (6.15) and the theorem 3.9, it follows that $a$ is left zero divisor. □

**Example 6.4.** The requirement that, for any $x$ be a root of the equation

$$axb = 0$$

$a$ must be a left zero divisor of $F$-algebra $A$, is necessary but not sufficient.

For instance, consider $R$-algebra of matrices $n \times n$. Consider matrices

$$E_k^i = (\delta_j^i \delta_k^l)$$

$$X = (x_k^i)$$

Then

$$E_j^i X E_l^k = (\delta_b^i \delta_j^a x_c^b \delta_d^k \delta_l^c) = (\delta_j^a x_l^i \delta_d^k) = x_l^i E_j^k$$

Therefore, the polynomial

$$p(X) = E_j^i X E_l^k$$

equals 0 iff $x_l^i = 0$.  □

Based on statements considered in this section we can assume that *homogeneous polynomial of degree 1 does not vanish identically.* However this statement requires more research.

It is evident that the theorem 5.9 is important. I recall that to prove this theorem we have assumed that the factors are not zero divisors.[14] However, if leading coefficient of polynomials are zero divisors, then conclusion of the theorem may not be true.

**Theorem 6.5.** *Let $A$-number $a$ be left zero divisor of $D$-algebra $A$. Let $p(x) \in A[x]$. Then the polynomial $p(x)a$ is left zero divisor of $A$-algebra $A[x]$.*

*Proof.* Since $p(x) \in A$, then the theorem follows from the theorem 3.6.  □

**Theorem 6.6.** *Let $a$ be right zero divisor of $D$-algebra $A$. Let $b$ be left zero divisor of $D$-algebra $A$. Let $ab \ne 0$. Then the polynomial $ap(x)b$ is left zero divisor of $A$-algebra $A[x]$.*

*Proof.* Since $p(x) \in A$, then the theorem follows from the theorem 3.3.  □

**Theorem 6.7.** *Let $a \in A^{\otimes 2}$ be nonsingular tensor. If we consider the equation (6.2) as transformation of algebra $A$, then we can write the inverse transformation in form*

$$(6.19) \qquad x = c^{pq} e_p b e_q$$

*where components $c^{pq}$ satisfy to equation*

$$(6.20) \qquad \delta_0^r \delta_0^s = a^{ij} c^{pq} C_{ip}^r C_{qj}^s$$

*Proof.* The theorem follows from the theorem [5]-4.2.  □

**Definition 6.8.** Let $a \in A^{\otimes 2}$ be nonsingular tensor. The tensor

$$a^{-1} = c^{pq} e_p \otimes e_q$$

is called **tensor inverse to tensor** $a$.                    $\square$

**Theorem 6.9.** *Let* $p(x) = p_1 \circ x$ *be homogeneous polynomial of power* $1$ *and* $p_1$ *be nonsingular tensor. Let*

$$r(x) = r_0 + r_1 \circ x + \ldots + r_k \circ x^k$$

*be polynomial of power* $k$*. Then*

$$\begin{aligned}
r(x) &= r_0 + q_{1\cdot0}p(x)q_{1\cdot1} + q_{2\cdot0}(x)p(x)q_{2\cdot1} + \ldots + q_{k\cdot0}(x)p(x)q_{k\cdot1} \\
&= r_0 + (q_{1\cdot0} \otimes q_{1\cdot1}) \circ p(x) + (q_{2\cdot0}(x) \otimes q_{2\cdot1}) \circ p(x) \\
&\quad + \ldots + (q_{k\cdot0}(x) \otimes q_{k\cdot1}) \circ p(x)
\end{aligned} \tag{6.21}$$

*Proof.* According to definitions 6.1, 6.8 and the theorem 6.7, the following equation is true

$$p_1^{-1} \circ p(x) = x \tag{6.22}$$

Based on the theorem 5.7, we can write the polynomial $r(x)$ as

$$\begin{aligned}
r(x) &= r_0 + r_{1\cdot0\cdot s}(xr_{1\cdot1\cdot s}) + (r_{2\cdot0\cdot s} \circ x)(xr_{2\cdot1\cdot s}) \\
&\quad + \ldots + (r_{k\cdot0\cdot s} \circ x^{k-1})(xr_{k\cdot1\cdot s}) \\
&= r_0 + r_{1\cdot0\cdot s}((1 \otimes r_{1\cdot1\cdot s}) \circ x) + (r_{2\cdot0\cdot s} \circ x)((1 \otimes r_{2\cdot1\cdot s}) \circ x) \\
&\quad + \ldots + (r_{k\cdot0\cdot s} \circ x^{k-1})((1 \otimes r_{k\cdot1\cdot s}) \circ x)
\end{aligned} \tag{6.23}$$

where

$$\begin{aligned}
r_1 &= r_{1\cdot0\cdot s}\underline{\otimes}(1 \otimes r_{1\cdot1\cdot s}) \quad r_{1\cdot0\cdot s} \in A \qquad r_{1\cdot1\cdot s} \in A \\
r_2 &= r_{2\cdot0\cdot s}\underline{\otimes}(1 \otimes r_{2\cdot1\cdot s}) \quad r_{2\cdot0\cdot s} \in A^{\otimes 2} \quad r_{2\cdot1\cdot s} \in A \\
&\ldots\ldots \\
r_k &= r_{k\cdot0\cdot s}\underline{\otimes}(1 \otimes r_{k\cdot1\cdot s}) \quad r_{k\cdot0\cdot s} \in A^{\otimes k} \quad r_{k\cdot1\cdot s} \in A
\end{aligned}$$

From the equations (6.22), (6.23), it follows that

$$\begin{aligned}
r(x) &= r_0 + r_{1\cdot0\cdot s}((1 \otimes r_{1\cdot1\cdot s}) \circ (p^{-1} \circ p(x))) \\
&\quad + (r_{2\cdot0\cdot s} \circ x)((1 \otimes r_{2\cdot1\cdot s}) \circ (p^{-1} \circ p(x))) \\
&\quad + \ldots + (r_{k\cdot0\cdot s} \circ x^{k-1})((1 \otimes r_{k\cdot1\cdot s}) \circ (p^{-1} \circ p(x))) \\
&= r_0 + r_{1\cdot0\cdot s}(((1 \otimes r_{1\cdot1\cdot s}) \circ p^{-1}) \circ p(x)) \\
&\quad + (r_{2\cdot0\cdot s} \circ x)(((1 \otimes r_{2\cdot1\cdot s}) \circ p^{-1}) \circ p(x)) \\
&\quad + \ldots + (r_{k\cdot0\cdot s} \circ x^{k-1})(((1 \otimes r_{k\cdot1\cdot s}) \circ p^{-1}) \circ p(x))
\end{aligned} \tag{6.24}$$

Let

$$p^{-1} = p'_{0\cdot t} \otimes p'_{1\cdot t} \tag{6.25}$$

Then

$$(1 \otimes r_{i\cdot1\cdot s}) \circ p^{-1} = (1 \otimes r_{i\cdot1\cdot s}) \circ (p'_{0\cdot t} \otimes p'_{1\cdot t}) = p'_{0\cdot t} \otimes p'_{1\cdot t}r_{i\cdot1\cdot s} \tag{6.26}$$

From the equations (6.24), (6.26), it follows that

$$
\begin{aligned}
r(x) &= r_0 + r_{1 \cdot 0 \cdot s}((p'_{0 \cdot t} \otimes p'_{1 \cdot t} r_{1 \cdot 1 \cdot s}) \circ p(x)) \\
&\quad + (r_{2 \cdot 0 \cdot s} \circ x)((p'_{0 \cdot t} \otimes p'_{1 \cdot t} r_{2 \cdot 1 \cdot s}) \circ p(x)) \\
&\quad + ... + (r_{k \cdot 0 \cdot s} \circ x^{k-1})((p'_{0 \cdot t} \otimes p'_{1 \cdot t} r_{k \cdot 1 \cdot s}) \circ p(x)) \\
&= r_0 + (r_{1 \cdot 0 \cdot s} \underline{\otimes} p'_{0 \cdot t})(p(x) p'_{1 \cdot t} r_{1 \cdot 1 \cdot s}) \\
&\quad + ((r_{2 \cdot 0 \cdot s} \underline{\otimes} p'_{0 \cdot t}) \circ x)(p(x) p'_{1 \cdot t} r_{2 \cdot 1 \cdot s}) \\
&\quad + ... + ((r_{k \cdot 0 \cdot s} \underline{\otimes} p'_{0 \cdot t}) \circ x^{k-1})(p(x) p'_{1 \cdot t} r_{k \cdot 1 \cdot s})
\end{aligned}
\tag{6.27}
$$

Let

$$
\begin{aligned}
q_{1 \cdot 0} &= r_{1 \cdot 0 \cdot s} \underline{\otimes} p'_{0 \cdot t} & q_{1 \cdot 1} &= p'_{1 \cdot t} r_{1 \cdot 1 \cdot s} \\
q_{2 \cdot 0}(x) &= (r_{2 \cdot 0 \cdot s} \underline{\otimes} p'_{0 \cdot t}) \circ x & q_{2 \cdot 1} &= p'_{1 \cdot t} r_{2 \cdot 1 \cdot s} \\
&\quad ...... & &\quad ...... \\
q_{k \cdot 0}(x) &= (r_{k \cdot 0 \cdot s} \underline{\otimes} p'_{0 \cdot t}) \circ x^{k-1} & q_{k \cdot 1} &= p'_{1 \cdot t} r_{k \cdot 1 \cdot s}
\end{aligned}
\tag{6.28}
$$

The equation (6.21) follows from equations (6.27), (6.28).                                    □

**Theorem 6.10.** *Let*

$$
p(x) = p_0 + p_1 \circ x
\tag{6.29}
$$

*be polynomial of power* 1 *and* $p_1$ *be nonsingular tensor. Let*

$$
r(x) = r_0 + r_1 \circ x + ... + r_k \circ x^k
$$

*be polynomial of power* $k$. *Then*[15]

$$
\begin{aligned}
r(x) &= r_0 - ((r_{1 \cdot 0 \cdot s} \otimes r_{1 \cdot 1 \cdot s}) \circ p_1^{-1}) \circ p_0 \\
&\quad - (((r_{2 \cdot 0 \cdot s} \circ x) \otimes r_{2 \cdot 1 \cdot s}) \circ p_1^{-1}) \circ p_0 \\
&\quad - ... - (((r_{k \cdot 0 \cdot s} \circ x^{k-1}) \otimes r_{k \cdot 1 \cdot s}) \circ p_1^{-1}) \circ p_0 \\
&\quad + ((r_{1 \cdot 0 \cdot s} \otimes r_{1 \cdot 1 \cdot s}) \circ p_1^{-1}) \circ p(x) \\
&\quad + (((r_{2 \cdot 0 \cdot s} \circ x) \otimes r_{2 \cdot 1 \cdot s}) \circ p_1^{-1}) \circ p(x) \\
&\quad + ... + (((r_{k \cdot 0 \cdot s} \circ x^{k-1}) \otimes r_{k \cdot 1 \cdot s}) \circ p_1^{-1}) \circ p(x) \\
&= r_0 - ((r_{1 \cdot 0 \cdot s} \otimes r_{1 \cdot 1 \cdot s} + (r_{2 \cdot 0 \cdot s} \circ x) \otimes r_{2 \cdot 1 \cdot s} \\
&\quad + ... + (r_{k \cdot 0 \cdot s} \circ x^{k-1}) \otimes r_{k \cdot 1 \cdot s}) \circ p_1^{-1}) \circ p_0 \\
&\quad + ((r_{1 \cdot 0 \cdot s} \otimes r_{1 \cdot 1 \cdot s} + (r_{2 \cdot 0 \cdot s} \circ x) \otimes r_{2 \cdot 1 \cdot s} \\
&\quad + ... + (r_{k \cdot 0 \cdot s} \circ x^{k-1}) \otimes r_{k \cdot 1 \cdot s}) \circ p_1^{-1}) \circ p(x)
\end{aligned}
\tag{6.30}
$$

*Proof.* From the equation (6.29) it follows that

$$
p_1 \circ x = -p_0 + p(x)
\tag{6.31}
$$

According to definitions 6.1, 6.8 and the theorem 6.7, from the equation (6.31) it is follows that

$$
p_1^{-1} \circ (-p_0 + p(x)) = x
\tag{6.32}
$$

---

[15]The style of the equation (6.30) is different from the style of the equation (6.21). I just want to show that we can use different styles to represent a polynomial.

Based on the theorem 5.7, we can write the polynomial $r(x)$ as

$$
\begin{aligned}
(6.33)\quad r(x) &= r_0 + r_{1\cdot 0\cdot s}(xr_{1\cdot 1\cdot s}) + (r_{2\cdot 0\cdot s} \circ x)(xr_{2\cdot 1\cdot s}) \\
&\quad + ... + (r_{k\cdot 0\cdot s} \circ x^{k-1})(xr_{k\cdot 1\cdot s}) \\
&= r_0 + (r_{1\cdot 0\cdot s} \otimes r_{1\cdot 1\cdot s}) \circ x + ((r_{2\cdot 0\cdot s} \circ x) \otimes r_{2\cdot 1\cdot s}) \circ x \\
&\quad + ... + ((r_{k\cdot 0\cdot s} \circ x^{k-1}) \otimes r_{k\cdot 1\cdot s}) \circ x
\end{aligned}
$$

where

$$
\begin{aligned}
r_1 &= r_{1\cdot 0\cdot s}\underline{\otimes}(1 \otimes r_{1\cdot 1\cdot s}) & r_{1\cdot 0\cdot s} &\in A & r_{1\cdot 1\cdot s} &\in A \\
r_2 &= r_{2\cdot 0\cdot s}\underline{\otimes}(1 \otimes r_{2\cdot 1\cdot s}) & r_{2\cdot 0\cdot s} &\in A^{\otimes 2} & r_{2\cdot 1\cdot s} &\in A \\
&... ... \\
r_k &= r_{k\cdot 0\cdot s}\underline{\otimes}(1 \otimes r_{k\cdot 1\cdot s}) & r_{k\cdot 0\cdot s} &\in A^{\otimes k} & r_{k\cdot 1\cdot s} &\in A
\end{aligned}
$$

From the equations (6.32), (6.33), it follows that

$$
\begin{aligned}
(6.34)\quad r(x) &= r_0 + (r_{1\cdot 0\cdot s} \otimes r_{1\cdot 1\cdot s}) \circ (p_1^{-1} \circ (-p_0 + p(x))) \\
&\quad + ((r_{2\cdot 0\cdot s} \circ x) \otimes r_{2\cdot 1\cdot s}) \circ (p_1^{-1} \circ (-p_0 + p(x))) \\
&\quad + ... + ((r_{k\cdot 0\cdot s} \circ x^{k-1}) \otimes r_{k\cdot 1\cdot s}) \circ (p_1^{-1} \circ (-p_0 + p(x))) \\
&= r_0 + ((r_{1\cdot 0\cdot s} \otimes r_{1\cdot 1\cdot s}) \circ p_1^{-1}) \circ (-p_0 + p(x)) \\
&\quad + (((r_{2\cdot 0\cdot s} \circ x) \otimes r_{2\cdot 1\cdot s}) \circ p_1^{-1}) \circ (-p_0 + p(x)) \\
&\quad + ... + (((r_{k\cdot 0\cdot s} \circ x^{k-1}) \otimes r_{k\cdot 1\cdot s}) \circ p_1^{-1}) \circ (-p_0 + p(x))
\end{aligned}
$$

The equation (6.30) follows from the equation (6.34).                                     □

The theorem 6.9 states that, for given homogeneus polynomial $p(x)$ of power 1 and given polynomial $r(x)$ of power $k$, we can represent the polynomial $r(x)$ as sum of products of the polynomial $p(x)$ over polynomials of power less than $k$. There is similar statement in the theorem 6.10 for given polynomial $p(x)$ of power 1.

Since the product in $D$-algebra $A$ is non commutative, then we can tell that the polynomial $p(x)$ is either left divisor of the polynomial $r(x)$, if

$$
r(x) = p(x)q(x)
$$

or right divisor of the polynomial $r(x)$, if

$$
r(x) = q(x)p(x)
$$

However, we can generalize this definition.

**Definition 6.11.** The polynomial $p(x)$ is called **divisor of polynomial** $r(x)$, if we can represent the polynomial $r(x)$ as

$$
(6.35)\qquad r(x) = q_{i\cdot 0}(x)p(x)q_{i\cdot 1}(x) = (q_{i\cdot 0}(x) \otimes q_{i\cdot 1}(x)) \circ p(x)
$$

□

## 8. Index

## 9. Special Symbols and Notations

# Многочлен над ассоциативной $D$-алгеброй

Александр Клейн

**Аннотация.** В статье рассмотрена алгебра полиномов над ассоциативной $D$-алгеброй с единицей. Использование тензорной записи позволяет упростить представление многочлена. Рассмотрены вопросы делимости многочлена произвольного порядка на многочлен первого порядка.

## Содержание

## 1. Предисловие

Теория многочленов над коммутативным кольцом (смотри определение многочлена, например, [1], с. 133) имеет утверждения, подобные утверждениям из теории чисел. К этим утверждениям относятся теорема о единственности разложения многочлена в произведение неприводимых многочленов ([3], с. 292), теорема о делении с остатком ([1], с. 141, теорема 2).

Теория многочленов над некоммутативной алгеброй неизмеримо сложнее ([11], p. 48). В этой статье мы сделали попытку немного продвинуться в этом направлении.

Возможность представления многочлена в виде

$$p(x) = a_0 + \sum_{k=1}^{n} a_k \circ x^k$$

является важным утверждением, из которого следуют многие утверждения статьи. Это утверждение основано на теореме 4.6 и её следствии 4.7.

Вначале концепция тензорного представления отображения свободной алгебры была применена к записи линейного отображения свободной алгебры над коммутативным кольцом.[1] Этот формат записи позволил сделать утверждения некоммутативного математического анализа более простыми.

Возможность представить многочлен с помощью тензоров позволяет избавиться от сложности и увидеть важные свойства многочлена. Опираясь на теоремы, рассмотренные в [5], я доказал теорему 6.9. Я надеюсь, это первый шаг к изучению делимости многочленов.

Александр Ложже был первым читателем моей статьи. Я благодарен ему за его ценные комментарии.

## 2. Соглашения

**Соглашение 2.1.** *В выражении вида*

$$a_{s \cdot 0} x a_{s \cdot 1}$$

*предполагается сумма по индексу* $s$. ☐

**Соглашение 2.2.** *Пусть $A$ - свободная алгебра с конечным или счётным базисом. При разложении элемента алгебры $A$ относительно базиса $\overline{\overline{e}}$ мы пользуемся одной и той же корневой буквой для обозначения этого элемента и его координат. В выражении $a^2$ не ясно - это компонента разложения элемента $a$ относительно базиса или это операция возведения в степень. Для облегчения чтения текста мы будем индекс элемента алгебры выделять цветом. Например,*

$$a = a^i e_i$$

☐

**Соглашение 2.3.** *Очень трудно провести границу между модулем и алгеброй. Тем более, иногда в процессе построения мы должны сперва доказать, что множество $A$ является модулем, а потом мы доказываем, что это множество является алгеброй. Поэтому для записи координат элемента модуля мы также будем пользоваться соглашением 2.2.* ☐

**Соглашение 2.4.** *Элемент $\Omega$-алгебры $A$ называется $A$-числом. Например, комплексное число также называется $C$-числом, а кватернион называется $H$-числом.* ☐

## 3. Делитель нуля ассоциативной $D$-алгебры

Пусть $D$ - коммутативное кольцо и $A$ - ассоциативная $D$-алгебра с единицей.

**Определение 3.1.** Пусть $a, b \in A$, $a \neq 0$, $b \neq 0$. Если $ab = 0$, то $a$ называется **левым делителем нуля**, а $b$ называется **правым делителем нуля**.[2] Если левый делитель нуля $a \in A$ является правым делителем нуля, то $a$ называется **делителем нуля**. ☐

**Теорема 3.2.** *Пусть $a, b \in A$, $a \neq 0$, $b \neq 0$. Из равенства $ab = 0$ не следует равенство $ba = 0$.*[3]

---

[1]Смотри, например, разделы [6]-1, [10]-1.

[2]Смотри также определение [2]-10.17.

[3]Доказательство теоремы основано на замечании в [2] после определения 10.17.

*Доказательство.* Пусть $A$ является алгеброй $3 \times 3$ матриц. Пусть

$$E_{12} = \begin{pmatrix} 0 & 1 & 0 \\ 0 & 0 & 0 \\ 0 & 0 & 0 \end{pmatrix} \qquad E_{23} = \begin{pmatrix} 0 & 0 & 0 \\ 0 & 0 & 1 \\ 0 & 0 & 0 \end{pmatrix}$$

Очевидно, что

$$E_{23}E_{12} = \begin{pmatrix} 0 & 0 & 0 \\ 0 & 0 & 0 \\ 0 & 0 & 0 \end{pmatrix}$$

Тем не менее

$$E_{12}E_{23} = \begin{pmatrix} 0 & 0 & 1 \\ 0 & 0 & 0 \\ 0 & 0 & 0 \end{pmatrix}$$

Нетрудно убедиться, что обе матрицы $(E_{12}, E_{23})$ являются делителями нуля. $\square$

**Теорема 3.3.** *Пусть $a$ является правым делителем нуля $D$-алгебры $A$. Пусть $b$ является левым делителем нуля $D$-алгебры $A$. Пусть $ab \neq 0$. Тогда для любого $d \in A$, $adb$ является делителем нуля $D$-алгебры $A$.*

*Доказательство.* Так как $a$ является правым делителем нуля $D$-алгебры $A$, то существует $c \neq 0$ такой, что $ca = 0$. Тогда

$$c(adb) = (ca)(db) = 0(db) = 0$$

Следовательно, $adb$ является правым делителем нуля $D$-алгебры $A$.

Так как $b$ является левым делителем нуля $D$-алгебры $A$, то существует $c \neq 0$ такой, что $bc = 0$. Тогда

$$(adb)c = (ad)(bc) = (ad)0 = 0$$

Следовательно, $adb$ является левым делителем нуля $D$-алгебры $A$. $\square$

**Теорема 3.4.** *Существует $D$-алгебра, в которой левый делитель нуля не является правым делителем нуля.*[4]

*Доказательство.* Пусть $A$ - свободное $R$-векторное пространство, имеющее базис Гамеля $\overline{\overline{e}}$.[5] Рассмотрим $R$-алгебру линейных отображений[6] $\mathcal{L}(R; A; A)$. Пусть отображение $f \in \mathcal{L}(R; A; A)$ определено равенством

$$(3.1) \qquad \begin{cases} f \circ e_i = e_{i-1} & i > 1 \\ f \circ e_1 = 0 \end{cases}$$

---

[4]Смотри также пример [2]-10.16.

[5]Базис Гамеля был рассмотрен в определении [9]-2.3.1.

[6]Пусть $f, g \in \mathcal{L}(R; A; A)$. Сумма отображений $f$ и $g$ определена равенством

$$(f + g) \circ x = f \circ x + g \circ x$$

Произведение отображений $f$ и $g$ определено равенством

$$(f \circ g) \circ x = f \circ (g \circ x)$$

Пусть отображение $g \in \mathcal{L}(R; A; A)$ определено равенством

(3.2) $$g \circ e_i = e_{i+1}$$

Пусть отображение $p \in \mathcal{L}(R; A; A)$ определено равенством

(3.3) $$\begin{cases} p \circ e_i = 0 & i > 1 \\ p \circ e_1 = e_1 \end{cases}$$

Из равенств (3.1), (3.2) следует, что[7]

(3.4) $$f \circ g = 1$$

3.4.1: Из равенств (3.1), (3.3) следует, что

$$(f \circ p) \circ e_1 = f \circ (p \circ e_1) = f \circ e_1 = 0$$
$$(f \circ p) \circ e_i = f \circ (p \circ e_i) = f \circ 0 = 0 \qquad i > 1$$

Следовательно, отображение $f$ является левым делителем нуля.

3.4.2: Пусть $h \in \mathcal{L}(R; A; A)$ - такое отображение, что

(3.5) $$h \circ f = 0$$

Из равенств (3.4), (3.5) следует, что

$$0 = 0 \circ g = (h \circ f) \circ g = h \circ (f \circ g) = h \circ 1 = h$$

Следовательно, отображение $f$ не является правым делителем нуля.

3.4.3: Из равенств (3.2), (3.3) следует, что

$$(p \circ g) \circ e_i = p \circ (g \circ e_i) = p \circ e_{i+1} = 0$$

Следовательно, отображение $g$ является правым делителем нуля.

3.4.4: Пусть $h \in \mathcal{L}(R; A; A)$ - такое отображение, что

(3.6) $$g \circ h = 0$$

Из равенств (3.4), (3.6) следует, что

$$0 = f \circ 0 = f \circ (g \circ h) = (f \circ g) \circ h = 1 \circ h = h$$

Следовательно, отображение $g$ не является левым делителем нуля.

$\square$

---

[7]В равенстве (3.4) мы видим, как меняются свойства матрицы, когда вместо матрицы с конечным числом строк и столбцов мы рассматриваем матрицу со счётным множеством строк и столбцов. Относительно базиса $\overline{\overline{e}}$, отображение $f$ имеет матрицу

$$f = \begin{pmatrix} 0 & 1 & 0 & 0 & ... \\ 0 & 0 & 1 & 0 & ... \\ 0 & 0 & 0 & 1 & ... \\ ... & ... & ... & ... & ... \end{pmatrix}$$

Первый столбец матрицы $f$ линейно зависит от остальных столбцов матрицы. Тем не менее, матрица $f$ обратима справа.

По-видимому, это явление связано с утверждением, что счётное множество имеет собственное счётное подмножество.

**Теорема 3.5.** *Существует $D$-алгебра, в которой левый делитель нуля является обратимым справа.*[8]

*Доказательство.* Пусть $\mathcal{L}(R; A; A)$ - $R$-алгебра, рассмотренная в доказательстве теоремы 3.4 Пусть отображение $f \in \mathcal{L}(R; A; A)$ определено равенством (3.1). Пусть отображение $g \in \mathcal{L}(R; A; A)$ определено равенством (3.2). Согласно утверждению 3.4.1, отображение $f$ является левым делителем нуля. Согласно равенству (3.4), отображение $f$ обратимо справа. $\square$

**Теорема 3.6.** *Пусть $a$ - левый делитель нуля $D$-алгебры $A$. Ненулевой элемент правого идеала $Aa$ является левым делителем нуля $D$-алгебры $A$.*

*Доказательство.* Согласно определению 3.1, существует $b \in A$, $b \neq 0$, такой, что $ab = 0$. Тогда для любого $c \in A$

$$(ca)b = c(ab) = c0 = 0$$

Следовательно, если $ca \neq 0$, то $ca$ является левым делителем нуля. $\square$

**Теорема 3.7.** *Если мы можем представить левый делитель нуля $a$ $D$-алгебры $A$ в виде произведения $a = cd$, то либо $c$, либо $d$ является левым делителем нуля.*

*Доказательство.* Если $d$ является левым делителем нуля, то, согласно теореме 3.6, $a$ является левым делителем нуля. Поэтому, чтобы доказать теорему, рассмотрим случай, когда $d$ не является левым делителем нуля. Согласно определению 3.1, существует $b \in A$, $b \neq 0$, такой, что $ab = 0$. Тогда

$$0 = ab = (cd)b = c(db)$$

Согласно определению 3.1, $db \neq 0$. Следовательно, $c$ является левым делителем нуля. $\square$

**Теорема 3.8.** *Пусть ни $a \in A$, ни $b \in A$ не являются левыми делителями нуля $D$-алгебры $A$. Тогда их произведение $ab$ не является левым делителем нуля $D$-алгебры $A$.*

*Доказательство.* Если произведение $ab$ является левым делителем нуля $D$-алгебры $A$, то, согласно теореме 3.7, либо $a$, либо $b$ является левым делителем нуля. Это противоречие доказывает теорему. $\square$

Чтобы лучше представить структуру множества делителей нуля $D$-алгебры $A$, мы рассмотрим следующую теорему.

**Теорема 3.9.** *Пусть $A$ - конечно мерная $D$-алгебра. Пусть $\overline{\overline{e}}$ - базис $D$-алгебры $A$. Пусть $C_{kl}^i$ - структурные константы $D$-алгебры $A$ относительно базиса $\overline{\overline{e}}$. Тогда координаты $a^i$ левого делителя нуля*

$$(3.7) \qquad a = a^i e_i$$

*удовлетворяют уравнению*

$$(3.8) \qquad \det \|C_{kl}^i a^k\| = 0$$

---

[8]Пусть $A$ - $D$-алгебра и $a \in A$ - левый делитель нуля. Тогда существует $b \in A$, $b \neq 0$, такой, что $ab = 0$. Пусть существует $c \in A$ такой, что $ca = 1$. Следовательно,

$$b = 1b = (ca)b = c(ab) = c0 = 0$$

Из полученного противоречия следует, что $a$ не является обратимым слева.

*Доказательство.* Если $A$-число $a \neq 0$ является левым делителем нуля, то согласно определению 3.1, $a \neq 0$ и существует $A$-число $b \neq 0$

$$b = b^i e_i$$

такое, что $ab = 0$. Следовательно, координаты $A$-чисел $a$ и $b$ удовлетворяют системе уравнений

$$(3.9) \qquad C_{kl}^i a^k b^l = 0$$

Если мы предположим, что $a$ известно, то мы можем рассматривать систему уравнений (3.9) как систему линейных уравнений относительно координат $b^l$. Уравнение (3.8) следует из утверждения, что число уравнений в системе линейных уравнений (3.9) равно числу неизвестных и система линейных уравнений (3.9) имеет нетривиальное решение. □

## 4. Многочлен над ассоциативной $D$-алгеброй

Пусть $D$ - коммутативное кольцо и $A$ - ассоциативная $D$-алгебра с единицей.

**Теорема 4.1.** *Пусть $p_k(x)$ - **одночлен степени** $k$ над $D$-алгеброй $A$.[9] Тогда*

4.1.1: *Одночлен степени 0 имеет вид* $p_0(x) = a_0$, $a_0 \in A$.

4.1.2: *Если $k > 0$, то*

$$p_k(x) = p_{k-1}(x) x a_k$$

*где $a_k \in A$.*

*Доказательство.* Мы докажем утверждение теоремы индукцией по степени $n$ одночлена.

Пусть $n = 0$. Так как одночлен $p_0(x)$ является константой, то мы получаем утверждение 4.1.1.

Пусть $n = k$. Последний множитель одночлена $p_k(x)$ является либо $a_k \in A$, либо имеет вид $x^l$, $l \geq 1$. В последнем случае мы положим $a_k = 1$. Множитель, предшествующий $a_k$, имеет вид $x^l$, $l \geq 1$. Мы можем представить этот множитель в виде $x^{l-1}x$. Следовательно, утверждение доказано. □

**Определение 4.2.** Обозначим $A_k[x]$ абелеву группу, порождённую множеством одночленов степени $k$. Элемент $p_k(x)$ абелевой группы $A_k[x]$ называется **однородным многочленом степени** $k$. □

*Замечание* 4.3. Согласно определению 4.2, однородный многочлен $p_k(x)$ степени $k$ является суммой одночленов степени $k$

$$p_k(x) = \sum_s p_{k \cdot s}(x)$$

Пусть

$$q_k(x) = \sum_t q_{k \cdot t}(x)$$

однородный многочлен степени $k$. Так как сложение в абелевой группе $A_k[x]$ коммутативно и ассоциативно, то сумма однородных многочленов $p$ и $q$ является суммой одночленов степени $k$

$$(p + q)(x) = p(x) + q(x) = \sum_s p_{k \cdot s}(x) + \sum_t q_{k \cdot t}(x)$$

---

[9]Аналогичное определение одночлена над телом рассмотрено в разделе [10]-16. Аналогичное определение одночлена над банаховой алгеброй рассмотрено в разделе [8]-5.2.

$\square$

**Теорема 4.4.** *Абелевая группа $A_k[x]$ является $A$-модулем.*

*Доказательство.* Пусть

$$p(x) = p_0 x p_1 ... p_{k-1} x p_k$$

одночлен. Для произвольного тензора $a \otimes b \in A \otimes A$, определено преобразование одночлена

(4.1)                    $(a \otimes b) \circ p(x) = ap(x)b = ap_0 x p_1 ... p_{k-1} x p_k b$

Множество преобразований (4.1) порождает представление $D$-алгебры $A \otimes A$ в абелевой группе $A_k[x]$. Следовательно, абелева группе $A_k[x]$ является $A$-модулем. $\square$

Множество одночленов степени $k$ не является базисом $A$-модуля $A_k[x]$. Например,

$$adxb = axdb \quad d \in D$$

Для многочлена, мы будем пользоваться записью

$$(a_0, a_1, ..., a_k) \circ x^k = a_0 x a_1 x ... x a_k$$

**Теорема 4.5.** *Отображение*

$$f : A^{k+1} \to A_k[x]$$

*определённое равенством*

$$f(a_0, a_1, ..., a_k) = (a_0, a_1, ..., a_k) \circ x^k$$

*является полилинейным отображением.*

*Доказательство.* Пусть $d \in D$. Из равенств

$$f(a_0, ..., da_i, ..., a_k) = a_0 x ... (da_i) ... x a_k = d(a_0 x ... a_i ... x a_k)$$
$$= df(a_0, ..., a_i, ..., a_k)$$
$$f(a_0, ..., a_i + b_i, ..., a_k) = a_0 x ... (a_i + b_i) ... x a_k$$
$$= (a_0 x ... a_i ... x a_k) + (a_0 x ... b_i ... x a_k)$$
$$= f(a_0, ..., a_i, ..., a_k) + f(a_0, ..., b_i, ..., a_k)$$

следует что отображение $f$ линейно по $a_i$. Следовательно, отображение $f$ - полилинейное отображение. $\square$

**Теорема 4.6.** *Определено линейное отображение*

$$(a_0 \otimes a_1 \otimes ... \otimes a_k) \circ x^k = (a_0, a_1, ..., a_k) \circ x^k$$

*Доказательство.* Утверждение теоремы является следствием теорем [7]-3.6.4, 4.5. $\square$

**Следствие 4.7.** *Однородный многочлен $p(x)$ может быть записан в виде*

$$p(x) = a_k \circ x^k \quad a_k \in A^{\otimes (k+1)}$$

$\square$

**Определение 4.8.** Обозначим

$$A[x] = \bigoplus_{n=0}^{\infty} A_n[x]$$

прямую сумму[10] $A$-модулей $A_n[x]$. Элемент $p(x)$ $A$-модуля $A[x]$ называется **многочленом** над $D$-алгеброй $A$.                                     □

Следовательно, многочлен степени $n$ может быть записан в виде

$$p(x) = a_0 + a_1 \circ x + ... + a_n \circ x^n \quad a_i \in A^{\otimes(i+1)} \quad i = 0, ..., n$$

**Определение 4.9.** Пусть

$$p(x) = a_0 + a_1 \circ x + ... + a_n \circ x^n \quad a_i \in A^{\otimes(i+1)} \quad i = 0, ..., n$$

многочлен степени $n$ над $D$-алгеброй $A$. $A^{\otimes(i+1)}$-число $a_i$ называется **коэффициентом многочлена** $p(x)$. $A^{\otimes(n+1)}$-число $a_n$ называется **старшим коэффициентом многочлена** $p(x)$.                                     □

## 5. Операции над многочленами

**Определение 5.1.** Пусть

$$p(x) = p_0 + p_1 \circ x + ... + p_n \circ x^n$$
$$q(x) = q_0 + q_1 \circ x + ... + q_n \circ x^n$$

многочлены степени $n$.[11] **Сумма многочленов** определена равенством

$$(p + q)(x) = p_0 + q_0 + (p_1 + q_1) \circ x + ... + (p_n + q_n) \circ x^n$$

□

**Определение 5.2.** Билинейное отображение

$$\underline{\otimes} \colon A^{\otimes n} \times A^{\otimes m} \to A^{\otimes(n+m-1)}$$

определено равенством

$$(a_1 \otimes ... \otimes a_n)\underline{\otimes}(b_1 \otimes ... \otimes b_n) = a_1 \otimes ... \otimes a_{n-1} \otimes a_n b_1 \otimes b_2 \otimes ... \otimes b_n$$

□

**Теорема 5.3.** *Однородный многочлен* $p \circ x^n$, $p \in A^{\otimes(n+1)}$, *порождает линейное отображение*

$$A_m[x] \to A_{n+m}[x]$$

*определённое равенством*

(5.1) $$(p \circ x^n) \circ (q \circ x^m) = (p\underline{\otimes}q) \circ x^{n+m}$$

---

[10] Смотри определение прямой суммы модулей в [1], страница 98. Согласно теореме 1 на той же странице, прямая сумма модулей существует.

[11] Пусть $q(x)$ - многочлен степени $m$, $m < n$,

$$q(x) = q_0 + q_1 \circ x + ... + q_m \circ x^m$$

Если мы положим $q_{m+1} = 0, ..., q_n = 0$, то мы можем рассматривать многочлен $q(x)$ как многочлен степени $n$.

*Доказательство.* Согласно следствию 4.7, $p \in A^{\otimes(n+1)}$, $q \in A^{\otimes(m+1)}$. Согласно определению 5.2,

$$(5.2) \qquad p\underline{\otimes}q \in A^{\otimes((n+1)+(m+1)-1)} = A^{\otimes((n+m)+1)}$$

Согласно следствию 4.7 и равенству (5.2), отображение (5.1) определено корректно. Согласно определению 5.2, из равенств

$$(p \circ x^n) \circ (r \circ x^m + s \circ x^m) = (p \circ x^n) \circ ((r + s) \circ x^m)$$
$$= (p\underline{\otimes}(r + s)) \circ x^{n+m}$$
$$= (p\underline{\otimes}r) \circ x^{n+m} + (p\underline{\otimes}s) \circ x^{n+m}$$
$$= (p \circ x^n) \circ (r \circ x^m) + (p \circ x^n) \circ (s \circ x^m)$$
$$(p \circ x^n) \circ (d(r \circ x^m)) = (p \circ x^n) \circ ((dr) \circ x^m)$$
$$= (p\underline{\otimes}(dr)) \circ x^{n+m} = d((p\underline{\otimes}r) \circ x^{n+m})$$
$$= d((p \circ x^n) \circ (r \circ x^m))$$

где $d \in D$, следует, что отображение (5.1) линейно. $\square$

*Замечание* 5.4. Пусть

$$p = p_0 \otimes p_1 \otimes ... \otimes p_n$$
$$r = r_0 \otimes r_1 \otimes ... \otimes r_m$$

Тогда равенство (5.1) можно записать следующим образом

$$(p_0x...xp_n) \circ (r_0x...xr_m) = ((p_0 \otimes ... \otimes p_n) \circ x^n) \circ ((r_0 \otimes ... \otimes r_m) \circ x^m)$$
$$= ((p_0 \otimes ... \otimes p_n)\underline{\otimes}(r_0 \otimes ... \otimes r_m)) \circ x^{n+m}$$
$$= (p_0 \otimes ... \otimes p_nr_0 \otimes ... \otimes r_m) \circ x^{n+m}$$
$$= p_0x...xp_nr_0x...xr_m$$

Следовательно, равенство (5.1) является определением произведения однородных многочленов. $\square$

**Определение 5.5.** Произведение однородных многочленов $p \circ x^n$, $r \circ x^m$ определено равенством

$$(p \circ x^n)(r \circ x^m) = (p\underline{\otimes}r) \circ x^{n+m}$$

$\square$

**Теорема 5.6.** *Пусть*

$$p(x) = a_k \circ x^k \qquad a_k \in A^{\otimes(k+1)}$$

*одночлен степени $k > 1$. Тогда многочлен $p(x)$ может быть представлен в одной из следующих форм*

$$(5.3) \qquad p(x) = (a_{k\cdot0} \circ x^{k-1})((1 \otimes a_{k\cdot1}) \circ x)$$

$$(5.4) \qquad p(x) = ((a_{k\cdot0} \circ x^{k-1}) \otimes a_{k\cdot1}) \circ x$$

*где*

$$(5.5) \qquad a_k = a_{k\cdot0}\underline{\otimes}(1 \otimes a_{k\cdot1}) \qquad a_{k\cdot0} \in A^{\otimes k} \qquad a_{k\cdot1} \in A$$

*Доказательство.* Опираясь на утверждение 4.1.2 и теорему 4.6, мы можем записать одночлен $p(x)$ в виде

$$(5.6) \qquad p_k(x) = (a_{k \cdot 0} \circ x^{k-1}) x a_{k \cdot 1}$$

Равенство (5.5) следует из определения 5.2 и равенства

$$a_k = a_{k \cdot 0} \otimes a_{k \cdot 1}$$

Так как для заданного значения $x$, выражение $a_{k \cdot 0} \circ x^{k-1}$ является $A$-числом, то равенство (5.4) следует из равенства (5.6) и теоремы 4.6. Равенство (5.3) следует из равенств (5.5), (5.6) и определения 5.5. $\qquad \square$

**Теорема 5.7.** *Пусть*

$$p(x) = a_k \circ x^k \quad a_k \in A^{\otimes(k+1)}$$

*однородный многочлен степени $k > 1$. Тогда многочлен $p(x)$ может быть представлен в одной из следующих форм*

$$(5.7) \qquad p(x) = (a_{k \cdot 0 \cdot s} \circ x^{k-1})((1 \otimes a_{k \cdot 1 \cdot s}) \circ x)$$

$$(5.8) \qquad p(x) = ((a_{k \cdot 0 \cdot s} \circ x^{k-1}) \otimes a_{k \cdot 1 \cdot s}) \circ x$$

*где*

$$(5.9) \qquad a_k = a_{k \cdot 0 \cdot s} \underline{\otimes}(1 \otimes a_{k \cdot 1 \cdot s}) \quad a_{k \cdot 0 \cdot s} \in A^{\otimes k} \quad a_{k \cdot 1 \cdot s} \in A$$

*Доказательство.* Согласно замечанию 4.3, однородный многочлен степени $k$ является суммой одночленов степени $k$. Теорема следует из теоремы 5.6, если рассмотреть индукцию по числу слагаемых. $\qquad \square$

*Замечание* 5.8. В этом разделе для нас несущественно, является ли многочлен $p(x)$ делителем нуля. Вопрос о делителях нуля $A$-алгебры $A[x]$ мы рассмотрим в разделе 6. При доказательстве теорем этого раздела, мы, не нарушая общности, будем предполагать, что рассматриваемые многочлены не являются делителями нуля. $\qquad \square$

Так как произведение однородных многочленов является билинейным отображением, то определение 5.5 может быть продолжено на произведение произвольных многочленов.

**Теорема 5.9.** *Пусть*

$$p(x) = p_0 + p_1 \circ x + ... + p_n \circ x^n$$
$$q(x) = q_0 + q_1 \circ x + ... + q_m \circ x^m$$

*многочлены. Если многочлен*

$$r(x) = r_0 + r_1 \circ x + ... + r_k \circ x^k$$

*является* **произведением многочленов**

$$(5.10) \qquad r(x) = p(x)q(x)$$

*то*

$$k = n + m$$

$$(5.11) \qquad r_h = \sum_{i+j=h} p_i \underline{\otimes} q_j \quad h = 0, ..., k \quad i \leq n \quad j \leq m$$

*Доказательство.* Мы докажем теорему индукцией по $n$, $m$.

- Следующая лемма является следствием определения 5.5.

*Лемма* 5.10. Теорема 5.9 верна для произведения однородных многочленов.

- Пусть теорема верна для однородного многочлена $p(x) = p_n \circ x^n$ степени $n$ и многочлена $q(x)$ степени $m = l$. Мы можем записать многочлен

$$q(x) = q_0 + q_1 \circ x + ... + q_l \circ x^l + q_{l+1} \circ x^{l+1}$$

степени $l + 1$ в виде

$$q(x) = q'(x) + q_{l+1} \circ x^{l+1}$$

где

$$q'(x) = q_0 + q_1 \circ x + ... + q_l \circ x^l$$

многочлен степени $l$. Так как произведение многочленов является билинейным отображением, то

$$(5.12) \qquad r(x) = (p_n \circ x^n)q(x) = (p_n \circ x^n)(q'(x) + q_{l+1} \circ x^{l+1})$$
$$= (p_n \circ x^n)q'(x) + (p_n \circ x^n)(q_{l+1} \circ x^{l+1})$$

Согласно предположению индукции, степень многочлена $(p_n \circ x^n)q'(x)$ равна $n + l$, а также

$$(5.13) \qquad r_h = \sum_{i+j=h} p_i \underline{\otimes} q_j \quad i = n \quad h = 0, ..., n + l$$

Согласно определению 5.5, степень многочлена $(p_n \circ x^n)(q_{l+1} \circ x^{l+1})$ равна

$$k = n + (l + 1)$$

а также

$$(5.14) \qquad r_h = \sum_{i+j=h} p_i \underline{\otimes} q_j \quad i = n \quad h = n + l + 1$$

Следовательно, утверждение теоремы верно для однородного многочлена $p(x) = p_n \circ x^n$ степени $n$ и многочлена $q(x)$ степени $m = l + 1$. Следовательно, мы доказали следующую лемму.

*Лемма* 5.11. Теорема 5.9 верна для произведения однородного многочлена $p(x)$ и многочлена $q(x)$.

Пусть теорема верна для многочлена $p(x)$ степени $n = l$ и многочлена $q(x)$ степени $m$. Мы можем записать многочлен

$$p(x) = p_0 + p_1 \circ x + ... + p_l \circ x^l + p_{l+1} \circ x^{l+1}$$

степени $l + 1$ в виде

$$p(x) = p'(x) + p_{l+1} \circ x^{l+1}$$

где

$$p'(x) = p_0 + p_1 \circ x + ... + p_l \circ x^l$$

многочлен степени $l$. Так как произведение многочленов является билинейным отображением, то

$$(5.15) \qquad \begin{aligned} r(x) &= p(x)q(x) = (p'(x) + p_{l+1} \circ x^{l+1})q(x) \\ &= p'(x)q(x) + (p_{l+1} \circ x^{l+1})q(x) \end{aligned}$$

Согласно предположению индукции, степень многочлена

$$r'(x) = p'(x)q(x)$$

равна $l + m$, а также

$$(5.16) \qquad r'_h = \sum_{i+j=h} p_i \underline{\otimes} q_j \quad i < l+1 \quad h = 0, ..., l+m$$

Согласно лемме 5.11, степень многочлена

$$r''(x) = (p_{l+1} \circ x^{l+1})q(x)$$

равна

$$k = (l+1) + m$$

а также

$$(5.17) \qquad r''_h = \sum_{i+j=h} p_i \underline{\otimes} q_j \quad i = l+1 \quad h = 0, ..., l+m+1$$

Так как

$$r(x) = r'(x) + r''(x)$$

то из определения 5.1 и из равенств (5.16), (5.17) следует, что

$$(5.18) \qquad r_h = \sum_{i+j=h} p_i \underline{\otimes} q_j \quad i \le l+1 \quad h = 0, ..., l+m+1$$

Следовательно, утверждение теоремы верно для многочлена $p(x)$ степени $n = l + 1$ и многочлена $q(x)$ степени $m$. $\qquad \square$

**Теорема 5.12.** *$A$-модуль* $A[x]$ *, оснащённый произведением* (5.10) *является $A$-алгеброй, называемой $A$-алгеброй многочленов над $D$-алгеброй $A$.*

*Доказательство.* Следствие определений [7]-3.2.1, 5.5 и теоремы 5.9. $\qquad \square$

## 6. Деление многочленов

В этом разделе мы будем предполагать, что алгебра $A$ определена над полем $F$ и является конечно мерной $F$-алгеброй.[12] Пусть $\overline{\overline{e}}$ - базис алгебры $A$ над полем $F$ и $C_{ij}^k$ - структурные константы алгебры $A$ относительно базиса $\overline{\overline{e}}$.

Рассмотрим линейное уравнение[13]

$$(6.1) \qquad a \circ x = b$$

---

[12]Из доказательства теоремы 3.4, мы видим, что утверждения линейной алгебры в векторном пространстве со счётным базисом отличаются от аналогичных утверждений в конечно мерном векторном пространстве. Поэтому необходимо дополнительное исследование, прежде чем мы можем сформулировать теоремы раздела 6 для $F$-алгебры со счётным базисом.

[13]Я рассматриваю решение уравнения (6.1) аналогично тому как я это сделал в разделе [5]-4.

где $a = a_{s\cdot 0} \otimes a_{s\cdot 1} \in A^{\otimes 2}$. Согласно теореме [4]-10.1.9 мы можем записать уравнение (6.1) в стандартной форме

$$a^{ij} e_i x e_j = b$$

(6.2)

$$a^{ij} = a_{s\cdot 0}{}^i_{\phantom{i}} a_{s\cdot 1}{}^j_{\phantom{j}} \quad a_{s\cdot 0} = a_{s\cdot 0}{}^i_{\phantom{i}} e_i \quad a_{s\cdot 1} = a_{s\cdot 1}{}^i_{\phantom{i}} e_i$$

Согласно теореме [4]-10.1.10 уравнение (6.2) эквивалентно уравнению

$$a^j_i x^i = b^j$$

(6.3)

$$a^j_i = a^{kr} C^p_{ki} C^j_{pr} \quad x = x^i e_i \quad b = b^i e_i$$

Согласно теории линейных уравнений над полем, если определитель

$$(6.4) \qquad \det \|a^j_i\| \neq 0$$

то уравнение (6.1) имеет единственное решение.

**Определение 6.1.** Тензор $a \in A^{\otimes 2}$ называется **невырожденным тензором**, если этот тензор удовлетворяет условию (6.4). □

**Теорема 6.2.** *Пусть $F$ - поле. Пусть $A$ - конечно мерная $F$-алгебра. Линейное уравнение*

$$(6.5) \qquad a \circ x = 0$$

*где $a = a_{s\cdot 0} \otimes a_{s\cdot 1} \in A^{\otimes 2}$. имеет корнем любое $x \in A$ тогда и только тогда, когда*

$$(6.6) \qquad a_{s\cdot 0}{}^k_{\phantom{k}} a_{s\cdot 1}{}^r_{\phantom{r}} C^p_{ki} C^j_{pr} = 0$$

*Доказательство.* Согласно теореме [4]-10.1.9 мы можем записать уравнение (6.5) в стандартной форме

$$a^{ij} e_i x e_j = 0$$

(6.7)

$$a^{ij} = a_{s\cdot 0}{}^i_{\phantom{i}} a_{s\cdot 1}{}^j_{\phantom{j}} \quad a_{s\cdot 0} = a_{s\cdot 0}{}^i_{\phantom{i}} e_i \quad a_{s\cdot 1} = a_{s\cdot 1}{}^i_{\phantom{i}} e_i$$

Согласно теореме [4]-10.1.10 уравнение (6.7) эквивалентно уравнению

$$a^j_i x^i = 0$$

(6.8)

$$a^j_i = a^{kr} C^p_{ki} C^j_{pr} \quad x = x^i e_i$$

Для того, чтобы любое $x \in A$ было корнем системы линейных уравнений (6.8), необходимо и достаточно, чтобы

$$(6.9) \qquad a^j_i = 0$$

Из равенств (6.8), (6.9) следует, что

$$(6.10) \qquad a^{kr} C^p_{ki} C^j_{pr} = 0$$

Равенство (6.6) является следствием равенств (6.7), (6.10). □

Очень трудно сказать, может ли условие (6.6) быть верным в какой-то алгебре. Однако мы можем проанализировать частный случай теоремы 6.2.

**Теорема 6.3.** *Пусть $F$ - поле. Пусть $A$ - конечно мерная $F$-алгебра. Для того, чтобы любое $x \in A$ было корнем уравнения*

$$axb = 0$$

*необходимо, чтобы $a$ являлся левым делителем нуля $F$-алгебры $A$.*

*Доказательство.* Пусть в уравнении (6.5) $s = 1$, $a_{1 \cdot 0} = a$, $a_{1 \cdot 1} = b$. Пусть $\overline{\overline{e}}$ - базис $F$-алгебры $A$. Тогда равенство (6.6) примет вид

$$(6.11) \qquad a^k b^r C_{ki}^p C_{pr}^j = 0$$

Если мы предположим, что $a$ известно, то мы можем рассматривать систему уравнений (6.11) как систему линейных уравнений относительно координат $b^r$. Число уравнений в системе линейных уравнений (6.11) равно числу неизвестных. Так как система линейных уравнений (6.11) имеет нетривиальное решение, то

$$(6.12) \qquad \det \|a^k C_{ki}^p C_{pr}^j\| = 0$$

для любого $r$. Из равенства (6.12) следует, что

$$(6.13) \qquad \det \|C_{pr}^j\| \det \|a^k C_{ki}^p\| = 0$$

Из равенства (6.13) следует, что либо

$$(6.14) \qquad \det \|C_{pr}^j\| = 0$$

либо

$$(6.15) \qquad \det \|a^k C_{ki}^p\| = 0$$

Пусть равенство (6.14) верно. Пусть $e_1 \in \overline{\overline{e}}$, $e_1 = 1$. Тогда

$$(6.16) \qquad e_p = e_p e_1 = C_{p1}^j e_j$$

Из равенства (6.14) следует, что существуют $c^p$, $c \neq 0$, такие, что

$$(6.17) \qquad c^p C_{p1}^j = 0$$

Из равенств (6.16), (6.17), следует, что

$$(6.18) \qquad c^p e_p = c^p C_{p1}^j e_j = 0$$

Следовательно, векторы $e_k$ линейно зависимы. Из противоречия следует, что в $F$-алгебре $A$ с единицей равенство (6.14) неверно.

Из равенства (6.15) и теоремы 3.9 следует, что $a$ является левым делителем нуля. $\qquad \square$

**Пример 6.4.** Требование, что $a$ является левым делителем нуля $F$-алгебры $A$, для того, чтобы любой $x$ был корнем уравнения

$$axb = 0$$

необходимо, но не достаточно.

Рассмотрим, например, $R$-алгебру матриц $n \times n$. Рассмотрим матрицы

$$E_k^i = (\delta_j^i \delta_k^l)$$

$$X = (x_k^i)$$

Тогда

$$E_j^i X E_l^k = (\delta_b^i \delta_j^a x_c^b \delta_d^k \delta_l^c) = (\delta_j^a x_l^i \delta_d^k) = x_l^i E_j^k$$

Следовательно, многочлен

$$p(X) = E^i_j X E^k_l$$

обращается в 0 тогда и только тогда, когда $x^i_l = 0$. □

Опираясь на утверждения, рассмотренные в этом разделе, можно предположить, что *однородный многочлен степени 1 не равен тождественно 0*. Однако это утверждение требует более детальный анализ.

Очевидно, что теорема 5.9 важна. Напомню, что при доказательстве этой теоремы мы предположили, что сомножители не являются делителями нуля.[14] Однако, если старшие коэффициенты многочлена являются делителями нуля, то заключение теоремы может оказаться неверным.

**Теорема 6.5.** *Пусть $A$-число $a$ является левым делителем нуля $D$-алгебры $A$. Пусть $p(x) \in A[x]$. Тогда многочлен $p(x)a$ является левым делителем нуля $A$-алгебры $A[x]$.*

*Доказательство.* Так как $p(x) \in A$, то теорема является следствием теоремы 3.6. □

**Теорема 6.6.** *Пусть $a$ является правым делителем нуля $D$-алгебры $A$. Пусть $b$ является левым делителем нуля $D$-алгебры $A$. Пусть $ab \neq 0$. Тогда многочлен $ap(x)b$ является левым делителем нуля $A$-алгебры $A[x]$.*

*Доказательство.* Так как $p(x) \in A$, то теорема является следствием теоремы 3.3. □

**Теорема 6.7.** *Пусть $a \in A^{\otimes 2}$ - невырожденный тензор. Если равенство (6.2) рассматривать как преобразование алгебры $A$, то обратное преобразование можно записать в виде*

$$(6.19) \qquad x = c^{pq} e_p b e_q$$

*где компоненты $c^{pq}$ удовлетворяют уравнению*

$$(6.20) \qquad \delta^r_0 \delta^s_0 = a^{ij} c^{pq} C^r_{ip} C^s_{qj}$$

*Доказательство.* Теорема является следствием теоремы [5]-4.2. □

**Определение 6.8.** Пусть $a \in A^{\otimes 2}$ - невырожденный тензор. Тензор

$$a^{-1} = c^{pq} e_p \otimes e_q$$

называется **тензор, обратный тензору** $a$. □

**Теорема 6.9.** *Пусть $p(x) = p_1 \circ x$ - однородный многочлен степени 1 и $p_1$ - невырожденный тензор. Пусть*

$$r(x) = r_0 + r_1 \circ x + ... + r_k \circ x^k$$

*многочлен степени $k$. Тогда*

$$r(x) = r_0 + q_{1 \cdot 0} p(x) q_{1 \cdot 1} + q_{2 \cdot 0}(x) p(x) q_{2 \cdot 1} + ... + q_{k \cdot 0}(x) p(x) q_{k \cdot 1}$$

$$(6.21) \qquad = r_0 + (q_{1 \cdot 0} \otimes q_{1 \cdot 1}) \circ p(x) + (q_{2 \cdot 0}(x) \otimes q_{2 \cdot 1}) \circ p(x)$$

$$+ ... + (q_{k \cdot 0}(x) \otimes q_{k \cdot 1}) \circ p(x)$$

---

[14]Смотри замечание 5.8.

*Доказательство.* Согласно определениям 6.1, 6.8 и теореме 6.7, следующее равенство верно

$$(6.22) \qquad p_1^{-1} \circ p(x) = x$$

Опираясь на теорему 5.7, мы можем записать многочлен $r(x)$ в виде

$$
\begin{aligned}
(6.23) \quad r(x) = {} & r_0 + r_{1 \cdot 0 \cdot s}(x r_{1 \cdot 1 \cdot s}) + (r_{2 \cdot 0 \cdot s} \circ x)(x r_{2 \cdot 1 \cdot s}) \\
& + ... + (r_{k \cdot 0 \cdot s} \circ x^{k-1})(x r_{k \cdot 1 \cdot s}) \\
= {} & r_0 + r_{1 \cdot 0 \cdot s}((1 \otimes r_{1 \cdot 1 \cdot s}) \circ x) + (r_{2 \cdot 0 \cdot s} \circ x)((1 \otimes r_{2 \cdot 1 \cdot s}) \circ x) \\
& + ... + (r_{k \cdot 0 \cdot s} \circ x^{k-1})((1 \otimes r_{k \cdot 1 \cdot s}) \circ x)
\end{aligned}
$$

где

$$
\begin{aligned}
r_1 &= r_{1 \cdot 0 \cdot s} \underline{\otimes} (1 \otimes r_{1 \cdot 1 \cdot s}) & r_{1 \cdot 0 \cdot s} &\in A & r_{1 \cdot 1 \cdot s} &\in A \\
r_2 &= r_{2 \cdot 0 \cdot s} \underline{\otimes} (1 \otimes r_{2 \cdot 1 \cdot s}) & r_{2 \cdot 0 \cdot s} &\in A^{\otimes 2} & r_{2 \cdot 1 \cdot s} &\in A \\
& ... ... \\
r_k &= r_{k \cdot 0 \cdot s} \underline{\otimes} (1 \otimes r_{k \cdot 1 \cdot s}) & r_{k \cdot 0 \cdot s} &\in A^{\otimes k} & r_{k \cdot 1 \cdot s} &\in A
\end{aligned}
$$

Из равенств (6.22), (6.23) следует, что

$$
\begin{aligned}
(6.24) \quad r(x) = {} & r_0 + r_{1 \cdot 0 \cdot s}((1 \otimes r_{1 \cdot 1 \cdot s}) \circ (p^{-1} \circ p(x))) \\
& + (r_{2 \cdot 0 \cdot s} \circ x)((1 \otimes r_{2 \cdot 1 \cdot s}) \circ (p^{-1} \circ p(x))) \\
& + ... + (r_{k \cdot 0 \cdot s} \circ x^{k-1})((1 \otimes r_{k \cdot 1 \cdot s}) \circ (p^{-1} \circ p(x))) \\
= {} & r_0 + r_{1 \cdot 0 \cdot s}(((1 \otimes r_{1 \cdot 1 \cdot s}) \circ p^{-1}) \circ p(x)) \\
& + (r_{2 \cdot 0 \cdot s} \circ x)(((1 \otimes r_{2 \cdot 1 \cdot s}) \circ p^{-1}) \circ p(x)) \\
& + ... + (r_{k \cdot 0 \cdot s} \circ x^{k-1})(((1 \otimes r_{k \cdot 1 \cdot s}) \circ p^{-1}) \circ p(x))
\end{aligned}
$$

Пусть

$$(6.25) \qquad p^{-1} = p'_{0 \cdot t} \otimes p'_{1 \cdot t}$$

Тогда

$$(6.26) \qquad (1 \otimes r_{i \cdot 1 \cdot s}) \circ p^{-1} = (1 \otimes r_{i \cdot 1 \cdot s}) \circ (p'_{0 \cdot t} \otimes p'_{1 \cdot t}) = p'_{0 \cdot t} \otimes p'_{1 \cdot t} r_{i \cdot 1 \cdot s}$$

Из равенств (6.24), (6.26) следует, что

$$
\begin{aligned}
(6.27) \quad r(x) = {} & r_0 + r_{1 \cdot 0 \cdot s}((p'_{0 \cdot t} \otimes p'_{1 \cdot t} r_{1 \cdot 1 \cdot s}) \circ p(x)) \\
& + (r_{2 \cdot 0 \cdot s} \circ x)((p'_{0 \cdot t} \otimes p'_{1 \cdot t} r_{2 \cdot 1 \cdot s}) \circ p(x)) \\
& + ... + (r_{k \cdot 0 \cdot s} \circ x^{k-1})((p'_{0 \cdot t} \otimes p'_{1 \cdot t} r_{k \cdot 1 \cdot s}) \circ p(x)) \\
= {} & r_0 + (r_{1 \cdot 0 \cdot s} \underline{\otimes} p'_{0 \cdot t})(p(x) p'_{1 \cdot t} r_{1 \cdot 1 \cdot s}) \\
& + ((r_{2 \cdot 0 \cdot s} \underline{\otimes} p'_{0 \cdot t}) \circ x)(p(x) p'_{1 \cdot t} r_{2 \cdot 1 \cdot s}) \\
& + ... + ((r_{k \cdot 0 \cdot s} \underline{\otimes} p'_{0 \cdot t}) \circ x^{k-1})(p(x) p'_{1 \cdot t} r_{k \cdot 1 \cdot s})
\end{aligned}
$$

Положим

$$(6.28) \quad \begin{aligned} q_{1\cdot 0} &= r_{1\cdot 0\cdot s}\underline{\otimes}p'_{0\cdot t} & q_{1\cdot 1} &= p'_{1\cdot t}r_{1\cdot 1\cdot s} \\ q_{2\cdot 0}(x) &= (r_{2\cdot 0\cdot s}\underline{\otimes}p'_{0\cdot t})\circ x & q_{2\cdot 1} &= p'_{1\cdot t}r_{2\cdot 1\cdot s} \\ &\dots\dots & &\dots\dots \\ q_{k\cdot 0}(x) &= (r_{k\cdot 0\cdot s}\underline{\otimes}p'_{0\cdot t})\circ x^{k-1} & q_{k\cdot 1} &= p'_{1\cdot t}r_{k\cdot 1\cdot s} \end{aligned}$$

Равенство (6.21) следует из равенств (6.27), (6.28). $\quad\square$

**Теорема 6.10.** *Пусть*

$$(6.29) \quad p(x) = p_0 + p_1 \circ x$$

*- многочлен степени 1 и $p_1$ - невырожденный тензор. Пусть*

$$r(x) = r_0 + r_1 \circ x + \dots + r_k \circ x^k$$

*многочлен степени $k$. Тогда*[15]

$$(6.30) \quad \begin{aligned} r(x) = {} & r_0 - ((r_{1\cdot 0\cdot s}\otimes r_{1\cdot 1\cdot s})\circ p_1^{-1})\circ p_0 \\ & - (((r_{2\cdot 0\cdot s}\circ x)\otimes r_{2\cdot 1\cdot s})\circ p_1^{-1})\circ p_0 \\ & - \dots - (((r_{k\cdot 0\cdot s}\circ x^{k-1})\otimes r_{k\cdot 1\cdot s})\circ p_1^{-1})\circ p_0 \\ & + ((r_{1\cdot 0\cdot s}\otimes r_{1\cdot 1\cdot s})\circ p_1^{-1})\circ p(x) \\ & + (((r_{2\cdot 0\cdot s}\circ x)\otimes r_{2\cdot 1\cdot s})\circ p_1^{-1})\circ p(x) \\ & + \dots + (((r_{k\cdot 0\cdot s}\circ x^{k-1})\otimes r_{k\cdot 1\cdot s})\circ p_1^{-1})\circ p(x) \\ = {} & r_0 - ((r_{1\cdot 0\cdot s}\otimes r_{1\cdot 1\cdot s} + (r_{2\cdot 0\cdot s}\circ x)\otimes r_{2\cdot 1\cdot s} \\ & + \dots + (r_{k\cdot 0\cdot s}\circ x^{k-1})\otimes r_{k\cdot 1\cdot s})\circ p_1^{-1})\circ p_0 \\ & + ((r_{1\cdot 0\cdot s}\otimes r_{1\cdot 1\cdot s} + (r_{2\cdot 0\cdot s}\circ x)\otimes r_{2\cdot 1\cdot s} \\ & + \dots + (r_{k\cdot 0\cdot s}\circ x^{k-1})\otimes r_{k\cdot 1\cdot s})\circ p_1^{-1})\circ p(x) \end{aligned}$$

*Доказательство.* Из равенства (6.29) следует, что

$$(6.31) \quad p_1 \circ x = -p_0 + p(x)$$

Согласно определениям 6.1, 6.8 и теореме 6.7, из равенства (6.31) следует

$$(6.32) \quad p_1^{-1}\circ(-p_0 + p(x)) = x$$

Опираясь на теорему 5.7, мы можем записать многочлен $r(x)$ в виде

$$(6.33) \quad \begin{aligned} r(x) = {} & r_0 + r_{1\cdot 0\cdot s}(xr_{1\cdot 1\cdot s}) + (r_{2\cdot 0\cdot s}\circ x)(xr_{2\cdot 1\cdot s}) \\ & + \dots + (r_{k\cdot 0\cdot s}\circ x^{k-1})(xr_{k\cdot 1\cdot s}) \\ = {} & r_0 + (r_{1\cdot 0\cdot s}\otimes r_{1\cdot 1\cdot s})\circ x + ((r_{2\cdot 0\cdot s}\circ x)\otimes r_{2\cdot 1\cdot s})\circ x \\ & + \dots + ((r_{k\cdot 0\cdot s}\circ x^{k-1})\otimes r_{k\cdot 1\cdot s})\circ x \end{aligned}$$

---

[15]Формат равенства (6.30) отличается от формата равенства (6.21). Я просто хочу показать, что мы можем использовать разные форматы для представления многочлена.

где

$$r_1 = r_{1\cdot 0\cdot s}\underline{\otimes}(1 \otimes r_{1\cdot 1\cdot s}) \quad r_{1\cdot 0\cdot s} \in A \qquad r_{1\cdot 1\cdot s} \in A$$

$$r_2 = r_{2\cdot 0\cdot s}\underline{\otimes}(1 \otimes r_{2\cdot 1\cdot s}) \quad r_{2\cdot 0\cdot s} \in A^{\otimes 2} \quad r_{2\cdot 1\cdot s} \in A$$

$$\dots\dots$$

$$r_k = r_{k\cdot 0\cdot s}\underline{\otimes}(1 \otimes r_{k\cdot 1\cdot s}) \quad r_{k\cdot 0\cdot s} \in A^{\otimes k} \quad r_{k\cdot 1\cdot s} \in A$$

Из равенств (6.32), (6.33) следует, что

$$
\begin{aligned}
(6.34) \quad r(x) = {}& r_0 + (r_{1\cdot 0\cdot s} \otimes r_{1\cdot 1\cdot s}) \circ (p_1^{-1} \circ (-p_0 + p(x))) \\
& + ((r_{2\cdot 0\cdot s} \circ x) \otimes r_{2\cdot 1\cdot s}) \circ (p_1^{-1} \circ (-p_0 + p(x))) \\
& + \dots + ((r_{k\cdot 0\cdot s} \circ x^{k-1}) \otimes r_{k\cdot 1\cdot s}) \circ (p_1^{-1} \circ (-p_0 + p(x))) \\
= {}& r_0 + ((r_{1\cdot 0\cdot s} \otimes r_{1\cdot 1\cdot s}) \circ p_1^{-1}) \circ (-p_0 + p(x)) \\
& + (((r_{2\cdot 0\cdot s} \circ x) \otimes r_{2\cdot 1\cdot s}) \circ p_1^{-1}) \circ (-p_0 + p(x)) \\
& + \dots + (((r_{k\cdot 0\cdot s} \circ x^{k-1}) \otimes r_{k\cdot 1\cdot s}) \circ p_1^{-1}) \circ (-p_0 + p(x))
\end{aligned}
$$

Равенство (6.30) следует из равенства (6.34). □

Теорема 6.9 утверждает, что, для заданного однородного многочлена $p(x)$ степени 1 и заданного многочлена $r(x)$ степени $k$, я могу представить многочлен $r(x)$ как сумму произведений многочлена $p(x)$ на многочлены степени меньше чем $k$. Аналогичное утверждение в теореме 6.10 для заданного многочлена $p(x)$ степени 1.

Поскольку произведение в $D$-алгебре $A$ некоммутативно, мы можем говорить, что многочлен $p(x)$ является либо левым делителем многочлена $r(x)$, если

$$r(x) = p(x)q(x)$$

либо правым делителем многочлена $r(x)$, если

$$r(x) = q(x)p(x)$$

Однако мы можем обобщить это определение.

**Определение 6.11.** Многочлен $p(x)$ называется **делителем многочлена** $r(x)$, если мы можем представить многочлен $r(x)$ в виде

$$(6.35) \qquad r(x) = q_{i\cdot 0}(x)p(x)q_{i\cdot 1}(x) = (q_{i\cdot 0}(x) \otimes q_{i\cdot 1}(x)) \circ p(x)$$

□

## 7. Список литературы

## 8. Предметный указатель

# 9. Специальные символы и обозначения